\documentclass[sn-mathphys,Numbered]{sn-jnl}

\usepackage{adjustbox}
\usepackage{geometry}
\usepackage{enumitem}

\usepackage{algorithm} 
\usepackage{algpseudocode} 
\usepackage{subcaption}
\usepackage{tikz}
\usetikzlibrary{shapes, arrows}
\usepackage{pgfplots}
\usepackage{amsmath}
\usetikzlibrary{patterns}
\usetikzlibrary{arrows}
\usepackage{lscape}
\usepackage{amssymb}
\usepackage{amsthm,amsmath,eurosym}

\usepackage{graphics}%
\usepackage{multirow}%
\usepackage{amsmath,amssymb,amsfonts}%
\usepackage{amsthm}%
\usepackage{mathrsfs}%
\usepackage[title]{appendix}%
\usepackage{xcolor}%
\usepackage{textcomp}%
\usepackage{manyfoot}%
\usepackage{booktabs}%
\usepackage{algorithm}%
\usepackage{algorithmicx}%
\usepackage{algpseudocode}%
\usepackage{listings}%
\usepackage{enumitem}
\usepackage{algorithm} 
\usepackage{algpseudocode} 
\usepackage{subcaption}
\usepackage{tikz}
\usetikzlibrary{shapes, arrows}
\usepackage{pgfplots}
\usepackage{amsmath}
\usetikzlibrary{patterns}
\usetikzlibrary{arrows}
\usepackage{lscape}
\usepackage{amssymb}
\usepackage{amsthm,amsmath,eurosym}

\theoremstyle{thmstyleone}%
\theoremstyle{thmstyletwo}%

\theoremstyle{thmstylethree}%

\begin{document}

\title[Article Title]{A hybrid meta-heuristic for the generation of feasible large-scale course timetables using instance decomposition}


\author[1]{\fnm{João} \sur{Almeida}}\email{joao.carvalho.almeida@tecnico.ulisboa.pt}

\author[1]{\fnm{José Rui} \sur{Figueira}}

\author[2]{\fnm{Alexandre P.} \sur{Francisco}}

\author[1]{\fnm{Daniel} \sur{Santos}}

\affil[1]{\orgdiv{CEGIST}, \orgname{Instituto Superior Técnico, Universidade de Lisboa} \orgaddress{\country{Portugal}}}

\affil[2]{\orgdiv{INESC-ID}, \orgname{Instituto Superior Técnico, Universidade de Lisboa} \orgaddress{\country{Portugal}}}


\abstract{This work presents a hybrid meta-heuristic for generating feasible course timetables, focusing on large-scale problems. We test it on instances from our university, for which the currently used commercial software takes several hours to find feasible solutions and often fails to comply with some constraints. The methodology includes elements from adaptive large neighbourhood search (ALNS), guided local search (GLS), variable neighbourhood search (VNS), and a novel instance decomposition technique. The constraint violations from different groups of constraints are regarded as objective functions to be minimised. The neighbourhood structures focus the search on the time slots where the most violations are identified. After a certain number of iterations without improvements, the most challenging constraint groups are given new weights, guiding the search towards non-dominated solutions, which improve the performance
of these objectives, even if the total sum of the effective violations increases. If this mechanism fails to escape these ``local optima'', a shaking phase is conducted. The decomposition mechanism operates as follows: curricula are iteratively introduced to the problem, and new feasible solutions are found considering the increasing set of lectures. The assignments from each iteration can be modified if needed in subsequent iterations. This methodology is tested in real-world instances from our university and in random sub-divisions of those. For sub-divisions with 400 curricula timetables, decomposition reduces the solution times to up to 27\%. For real-world instances, with 1288 curricula timetables, the reduction is 18\%. The experiments also reveal that clustering curricula with more lectures and professors in common when incrementing the instances improved the solution times by 18\% more than a random increment order. Using our methodology, feasible solutions for the real-world instance are found in 21 minutes, on average, whereas the commercial software takes several hours.}




\maketitle

\section{Introduction}
\noindent Course timetabling is one of the most studied educational scheduling problems \citep{TanEtAl2021}, { which also include exam scheduling \citep{LEITE2019} or thesis defence scheduling \citep{ALMEIDA2023}}. Moreover, it is an arduous and time-consuming burden for those tasked with its resolution at the university administration level \citep{BabaeiEtAl2015}. The problem can be defined through a 5W framework as follows: we want to schedule What lecture, taught by Who (professors), to Whom (curricula), Where (room) and When (day and hour).

A feasible timetable fulfils a particular set of constraints specific to each case \citep{ThepphakornPongchaoren2020}, such as students, professors, and rooms cannot have more than one lecture scheduled at the same time, some lectures can only be scheduled after another has taken place, or some time-slots are unavailable for certain professors, curricula or rooms \citep{ChenEtAl2021}.

Course timetabling impacts several stakeholders. On the one hand, there are those responsible for the scheduling process who have to do extra corrective work in short amounts of time if the timetabling software does not find satisfactory solutions in an acceptable time frame. On the other hand, there are stakeholders whose daily activities can be affected by the quality of the timetables. They can be students, professors, or managerial stakeholders, who must manage costs or room availability for other purposes.

Analysing real-world course timetabling systems shows a variety of concerns and points of view depending on the characteristics of each case. \cite{VermuytenEtAl2016} show a case where a campus reorganization led to congestion in the hallways and stairs, frequently causing lectures to start late. To solve this problem, the authors propose a course timetabling solution which minimises the travel times between consecutive classes for students. \cite{ThepphakornPongchaoren2020} provide a different perspective. Their solution approach minimises operational costs by rearranging the timetables to reduce, for example, professor hiring costs or cleaning and maintenance costs. \cite{BouffletEtAl2021} contribute with another point of view. Their previous timetabling methodology often failed to schedule all the lectures, leading to extraordinary managerial work to find ways to schedule the remaining lectures. Thus, they propose a solution that aims to maximise the number of scheduled lectures to reduce the burden on the management.

Our research is inspired by the course timetabling system of the Instituto Superior Técnico of the University of Lisbon (IST), where more than 4000 lectures must be scheduled each semester in two different campi. {Comparatively, the largest instances solved by the previously mentioned works have at most 1000 lectures}. Moreover, many of the lectures in our instances are shared between different curricula and classes. Additionally, the number of certain types of rooms is considerably limited for the number of lectures that must be scheduled. Thus, it is impossible to schedule the lectures of different departments separately as they compete for the same rooms and are taught simultaneously to multiple groups of students from different courses. Furthermore, IST operates using both semester and trimester academic structures, which poses a unique scheduling challenge. Certain courses span a semester, while others are taught in a single trimester.

Due to the size and intricacies of the problem, the commercial timetabling software currently used by IST cannot efficiently find timetables that comply with all the rules. For example, it is not uncommon that, in this initial schedule, a professor ends up with two lectures in different campi without any time in between or that some students cannot enrol in some courses due to juxtaposition with others. These problems must then be resolved by hand. {This burdensome process takes several weeks of manually rearranging individual timetables, checking their feasibility with the respective stakeholders, and booking new rooms.}

{This work proposes a new hybrid meta-heuristic for generating feasible course timetables. It integrates elements from adaptive large neighbourhood search (ALNS), guided local search (GLS), and variable neighbourhood search (VNS). Its objective function regards the minimisation of hard constraint violations. In our experiments, all instances were known to be feasible. Thus, no time limit was set, and the meta-heuristic stops when the value of the objective function reaches 0. If the feasibility of the instances was not known, a time or iteration limit could be set, and the meta-heuristic would still minimise the number of violations. The search is directed toward lectures and time slots with the highest number of violations using novel neighbourhood structures. When there is no improvement after a certain number of iterations, the most challenging constraint groups are assigned new weights. This guides the search towards non-dominated solutions that enhance the performance of these objectives, even if the total sum of violations increases. If this approach fails to move away from these ``local optima'' a shaking phase is implemented.}

{This work contributes to the literature by proposing a new hybrid meta-heuristic for generating feasible solutions for course timetabling problems tailored to solve large and complex instances}. For this purpose, we introduce novel neighbourhood structures for ALNS and features for GLS. The neighbourhood structures target specific lectures which are potentially causing violations. This is exceptionally useful when dealing with large and complex instances where random searches would most likely miss the problematic areas.
We also explore the use of new instance decomposition techniques in generating feasible course timetables.
Our decomposition mechanism functions as follows: a set of lectures belonging to specific curricula is added to the problem, and the meta-heuristic is used to find a new feasible solution considering the current set of lectures. When we have a feasible solution, we add the next set. Let us note that the assignments of each iteration can be changed if necessary in the following iterations. This mechanism allows us to simplify the problem by dealing with less complex states, with fewer violations, in each iteration.

The computational experiments are conducted on real-world instances and random sub-divisions of those. The decomposition mechanism promoted solution time reductions of up to 27\% compared to solving the problems without decomposition mechanisms in the sub-divisions and of up to 18\% in the real-world instance. The results also imply a significant positive impact if the curricula are added in clusters, \textit{i.e.}, curricula which share classes or professors are included in the problem at similar times. 

The remainder of this paper is structured in the following manner: in Section \ref{lr}, we provide a state-of-the-art review of recent developments in course timetabling problems.  In Section \ref{mpm}, we define the course timetabling problem for IST. In Section \ref{af}, we present the algorithmic framework of our meta-heuristic. In Section \ref{ce}, we analyse the computational experiments, including the IST instances and adaptions of the 2007 International Timetabling Competition curriculum-based course timetabling instances. Finally, in Section \ref{c}, we provide a conclusion and propose several new research paths.

\section{Literature Review}\label{lr}
\noindent In this section, we present a review of the latest developments in course timetabling. In Subsection \ref{lr-pd}, we summarise common problem definitions based on competition tracks and real-world problems. In Subsection \ref{lr-sm}, we present the solution approaches employed to solve such problems. In Subsection \ref{lr - fs}, we analyse other works proposing different methodologies for finding feasible course timetabling solutions.

\subsection{Problem definitions}\label{lr-pd}

\noindent The course timetabling problem has some characteristic constraints that, even if differently handled, are globally present in the literature \citep{ThepphakornPongchaoren2020, Herres2021}. These general constraints can be divided into two groups. The first group prevents student, professor, and room assignment juxtaposition, and the second group ensures that each lecture is assigned a time slot and room. Moreover, some problems follow other university regulations, leading to the addition of some problem-specific constraints. Contrariwise, objectives tend to be more varied between problems.

Standardised benchmark instances for the course timetabling problem have been widely studied \citep{BabaeiEtAl2015}. The two most studied benchmark instances were introduced in the 2007 International Timetabling Competition. They are the curriculum-based course timetabling problem \citep{Soria-Alcaraz2016, Kiefer2017,BaggerEtAl2019B,BaggerEtAl2019, PillayOzcan2019, LindahlEtAl2019, GülcüAkkan2020, AkkanEtAl2021, Rosa-Rivera2021, Song2021} and the post-enrolment course timetabling problem \citep{GohEtAl2017,Nagata2018, GohEtAl2019}. However, the 2019 International Timetabling Competition \citep{Sylejmani2022} introduced a new set of benchmark instances.

The curriculum-based course timetabling problem considers the general constraints and an additional constraint, which states that some curricula and professors are unavailable for certain time slots \citep{BettinelliEtAl2015}. Additionally, there are four objectives. Specifically, the adherence to a minimum number of working days for the lectures of a course, the promotion of compact schedules for each curriculum, respecting room capacity, and the minimisation of the number of rooms assigned to the lectures of the same course \citep{BettinelliEtAl2015}. Moreover, some works have also considered the robustness of the solutions as an additional objective \citep{AkkanGulcu2018, LindahlEtAl2019, GülcüAkkan2020, AkkanEtAl2021}.

Besides the general constraints, the post-enrolment course timetabling problem also considers room capacity violations, lecture precedence, and unavailable time slots as constraints \citep{Nagata2018}. Moreover, three objectives are considered. Specifically, avoiding specific time slots, limiting the number of consecutive lectures for a student, and avoiding days with a single lecture for a student \citep{Nagata2018}. 

The major difference between the problem of ITC2019 is the required student sectioning \citep{Sylejmani2022}. In this problem, the students or curricula are not already assigned to each lecture, and this assignment is one of the decisions the optimisation algorithms must make. Besides this concern, this problem also includes other constraints, such as maximum consecutive workloads, classes that must be taught on different days, or classes that must precede other classes.

{Real-world problems are also addressed in the literature. \cite{ThepphakornPongchaoren2020} deal with the course timetabling problem at the Faculty of Engineering of Naresuan University, and their largest instances have 1009 lectures and 66 curricula. Besides general constraints, this problem also considers constraints on certain room requirements, professor and student time slot availability, and consecutive lecture assignments. Their objective function aims to minimise operating costs for the faculty. Specifically, its components are the operating costs generated by assigning certain rooms, hiring lecturers for certain time slots, and cleaning and setting up rooms, which should be assigned as compact schedules as possible. \cite{BouffletEtAl2021} addresses the University of Technology of Compiègne course timetabling problem. Due to the size of their instances, they adopted the approach of including the assignment of lectures as an objective while maintaining the remaining general constraints. Their largest instances include 652 lectures.
Additionally, they also include avoiding building changes and minimising the number of necessary time slots. \cite{Wu2011} use the problem of the Department of Computer Science and Information Engineering at the National Changhua University of Education as their case study, solving instances with six curricula and 45-50 courses. Besides the general constraints, they also include constraints regarding room requirements and certain time-slot availability concerns. As objectives, they consider some preferred time slots, maximum consecutive and daily assignments, and avoidance of similar time slots for specific optional lectures. \cite{Shiau2011} solve a problem from a Taiwanese university, with at most 28 curricula. They include general and room requirement constraints, and their objectives include time slot preferences, course time slot distribution preferences, and minimisation of room preferences. \cite{Aladag2009} address instances from the Department of Statistics at Hacettepe University, with four curricula and 57 lectures. Their problem-specific constraints are room requirements and lectures scheduled on different days. Their objective functions include compactness and workload considerations, time slot preferences, and room capacity violations.}

\subsection{Solution improvement approaches}\label{lr-sm}

\noindent In course timetabling literature, solution approaches that provide a single solution \citep{Kiefer2017, Nagata2018, Bagger2018EtAl, BaggerEtAl2019B, PillayOzcan2019, GohEtAl2019, BaggerEtAl2019, ThepphakornPongchaoren2020, BouffletEtAl2021, AkkanEtAl2021, Rosa-Rivera2021, Song2021} are more often addressed than multi-objective methodologies \citep{AkkanGulcu2018, LindahlEtAl2018, LindahlEtAl2019, GülcüAkkan2020}, which aim to provide better knowledge about the different possible timetabling options. 

For finding a single solution, authors have studied meta-heuristic-based approaches \citep{Kiefer2017, AkkanGulcu2018, Nagata2018, GohEtAl2019, ThepphakornPongchaoren2020, AkkanEtAl2021, Rosa-Rivera2021, Song2021}, with the latter also including machine learning in their approach, mixed integer linear programming based approaches \citep{Bagger2018EtAl, BaggerEtAl2019, BouffletEtAl2021, BaggerEtAl2019B, Dunke2023}, and hyper-heuristics based approaches \citep{PillayOzcan2019}.

Some authors use multi-objective approaches \citep{AkkanGulcu2018, LindahlEtAl2018, LindahlEtAl2019, GülcüAkkan2020}. Moreover, all of these address the curriculum-based course timetabling problem and implement bi-criteria models. Works that introduced robustness metrics used this methodology to understand the trade-off between the weighted objective function and their respective robustness metrics \citep{AkkanGulcu2018, LindahlEtAl2019, GülcüAkkan2020}. Conversely, \cite{LindahlEtAl2018} use bi-criteria optimisation to assess trade-offs between pairs of objective functions already present in the original problem. Employed algorithms were based on meta-heuristics \citep{AkkanGulcu2018, GülcüAkkan2020} and mixed integer linear programming \citep{LindahlEtAl2018, LindahlEtAl2019}.

\subsection{Finding feasible solutions}\label{lr - fs}
\noindent {Most meta-heuristics consider initial solutions to be optimised as their starting points. The algorithms employed to generate these initial solutions can significantly impact the quality of the final solution. This influence can come from both the initial performance and the time required to generate the solutions.
}

{Several works have studied the generation of initial feasible solutions in different settings. \cite{Sousa2016} address the scheduling of energy resources in smart grids and compare the use of different strategies, such as random solutions, ant colony optimisation, naive-scheduling heuristics, pre-scheduling heuristics, and mixed integer linear programming. \cite{Juman2015} and \cite{Amaliah2022} propose heuristic approaches for transportation problems. \cite{Viana2022} uses genetic algorithms to generate feasible and diverse dungeon levels. }

{Regarding course timetabling, \cite{Song2018} proposes an iterated local search algorithm and tackles a problem considering student and room clashes, room requirements, and complete assignment of events as constraints. Their initial solution mechanism attempts to schedule as many lectures as possible without generating student and room clashes and respecting room requirements. They then use a simulated annealing procedure to improve this initial solution. }

{\cite{Pillay2019} proposes a hyper-heuristic approach.  They explore the automatic induction of two construction heuristics types: arithmetic and hierarchical heuristics. Genetic programming is employed to evolve arithmetic heuristics, while a combination of genetic programming, genetic algorithms, and the generation of random heuristic combinations is used to generate hierarchical heuristics.}

{\cite{Goh2019} proposes a Monte Carlo tree search methodology. Their problems include similar constraints to the ones present in \cite{Song2018}, with the addition of predefined time slots for some lectures and lecture ordering constraints. Besides their simulation approach, they provide several heuristic methods to improve the simulation efficiency by pruning their Monte Carlo trees.}

{Like \cite{Song2018}, our work uses a local search based approach but with a few key distinctions. Our initial solution approach schedules all the lectures while respecting student/curriculum-related hard constraints and disregarding professor, room, and lecture ordering constraints. Moreover, our neighbourhood structures specifically focus the search on the lectures which are causing the violations. Their objective function represents the number of violations. In our work, we include weights for different types of constraint violations. These weights are constantly updated during the search to penalise violations identified as the most challenging to overcome. Our largest tested instances include more than 4000 lectures, whereas \cite{Song2018} consider at most 1000. For instances of this dimension, instance decomposition and increment significantly improved solution times. As far as we know, this procedure was never attempted before in course timetabling problems.
}

\section{The course timetabling problem for Instituto Superior Técnico} \label{mpm}

\noindent {In this section, we define the course timetabling problem for IST. In Subsection \ref{indices}, we introduce the notation for the indices. In Subsection \ref{mpm p}, we provide the required parameters. In Subsection \ref{mpm dv}, we characterise the decision and auxiliary variables. In Subsection \ref{mpm c}, we explore the constraints that determine the feasibility of the timetables. In Subsection \ref{instances}, we present the instances solved in our work.}

\subsection{Indices}\label{indices}

\noindent The indices are presented in this subsection.

\begin{itemize}[label={--}]
    \item $i = 1,\ldots, n_i$, are the indices related to the lectures;
    \item $j = 1,\ldots, n_j$, are the indices related to the days;
    \item $k = 1,\ldots, n_k$, are the indices related to the hour slots in each day;
    \item $t = 1,\ldots, n_t$, are the indices related to the room types;
   \item $\ell = 1,\ldots, n_\ell$, are the indices related to the curricula;
    \item $p = 1,\ldots, n_p$, are the indices related to the professors;

    \item $\mu = 1,\ldots, n_\mu$, are the indices related to the campi;
\end{itemize}
\subsection{Parameters}\label{mpm p}
\noindent The parameters are presented in this subsection. They are divided into four groups. The first is related to associations between lectures, curricula, and professors. The second is related to lecture characteristics. The third is related to curricula and professor workloads.  The fourth is related to rooms. 
\begin{enumerate}
    \item \textit{Related to associations between lectures, curricula and professors.}
    \begin{itemize}[label={--}]
       \item { $a_{\ell} \in \{ 1,\ldots, n_i\}$, is the set of lectures, $i$, which is part of a curricula, $\ell$, for all $\ell = 1,\ldots, n_\ell$.}
      \item { $b_{p} \in \{ 1,\ldots, n_p\}$, is the set of lectures, $i$, taught by a professor, $p$, for all $p = 1,\ldots, n_p$.}
    \end{itemize}

    \item \textit{Related to lecture characteristics.}
    \begin{itemize}[label={--}]
        \item $c_{i} \in \mathbb{N}$, is the duration of a lecture, $i$, in hour slots, for all $i = 1,\ldots, n_i$.
        \item { $d_{i} \in \{ 1,\ldots, n_i\}$, is the set of lectures, $\overline{i}\neq i$, which cannot be taught on the same day as another lecture, $i$, for all $i = 1,\ldots, n_i$.}
        \item { $e_{i} \in \{ 1,\ldots, n_i\}$, is the set of lectures, $\overline{i}\neq i$, which must be taught before another lecture, $i$ for all $i = 1,\ldots, n_i$.}

    \end{itemize}

    \item \textit{Related to curricula and professor workloads.}
     \begin{itemize}[label={--}]
        \item $m_{\ell} \in \mathbb{N}$, is the maximum daily workload for a given curricula, $\ell$, in number of hour slots, for all $\ell = 1,\ldots, n_\ell$.
        \item $o_{p} \in \mathbb{N}$, is the maximum daily workload for a given professor, $p$, in number of hour slots, for all $p = 1,\ldots, n_p$.
        \item $r_{\ell} \in \mathbb{N}$, is the maximum consecutive workload for a given curricula, $\ell$, in number of hour slots, for all $\ell = 1,\ldots, n_\ell$.
        \item $s_{p} \in \mathbb{N}$, is the maximum consecutive workload for a given professor, $p$, in number of hour slots, for all $p = 1,\ldots, n_p$.
    \end{itemize}

    \item \textit{Related to campi.}
      \begin{itemize}[label={--}]
\item $q_t \in \mathbb{N}$, is the number of rooms of a type, $t$.
\item $v_{it} \in \{0,1\}$, is 1 if a lecture, $i$, is taught in a room of a type, $t$.
\item {$u_{\mu\overline{\mu}}\in \mathbb{N}$, is the number of time-slots necessary to travel between two campi, $\mu$ and $\overline{\mu}$.}
\item $g_{i} \in \{1,\ldots,n_\mu\}$, is the campi, $\mu$, in which a lecture, $i$, is taught.
      \end{itemize}
       
\end{enumerate}
\subsection{Variables}\label{mpm dv}
\noindent {The variables are presented in this subsection.}

\begin{enumerate}
    \item \textit{Decision variables.}
      \begin{itemize}[label={--}]
      \item $x_{ijk} \in \{0,1\}$, is 1 if a lecture, $i$, is taught in a day, $j$, and hour slot, $k$; and 0 otherwise, for all $i = 1,\ldots, n_i$, $j = 1,\ldots, n_j$, $k = 1,\ldots, n_k$.
      \end{itemize}

      \item \textit{Auxiliary variables.}
      \begin{itemize}[label={--}]
      \item $y_{jk\ell}\in \{0,1\}$, is 1 if a curriculum, $\ell$, has a lecture scheduled in a time slot, ($j,k$); and 0 otherwise, for all $j = 1,\ldots, n_j$, $k = 1,\ldots, n_k$, $\ell = 1,\ldots,n_\ell$.
    
    \item $\overline{y}_{jkp}\in \{0,1\}$, is 1 if a professor, $p$, has a lecture scheduled in a time slot, ($j,k$); and 0 otherwise, for all $j = 1,\ldots, n_j$, $k = 1,\ldots, n_k$, $p = 1,\ldots,n_p$.

 \end{itemize}
\end{enumerate}


\subsection{Constraints}\label{mpm c}
\noindent {The constraints that define the course timetabling problem at IST are presented in this subsection. They are divided into four groups. The first is related to curriculum, professor, and room availability. The second is related to curriculum and professor workloads. The third is related to campus changes. Finally, the fourth is related to lecture assignments.}

\begin{enumerate}
  
    \item \textit{Curriculum, professor and room availability}. These constraints ensure that students, professors and rooms are available at the time of their assignments.
    \begin{enumerate}
        \item \textit{Professor and student juxtaposition.} Students in a certain curriculum, $\ell$, or professors, $p$, do not have more than one lecture assigned at any given time.

        \item \textit{Room juxtaposition. } There are enough rooms of a type, $t$, for the lectures, $i$, which must be taught in a room of that type, $v_{it}=1$.

    \end{enumerate}
    
    \item \textit{Curriculum and professor workload}. These constraints ensure that regulations regarding the maximum student and professor workloads are met.
    \begin{enumerate}
        \item \textit{Maximum daily workload.} In a day, $j$, a curriculum, $\ell$, or professor, $p$, does not have more hour slots, $k$, of lectures, $i$, assigned than the regulations allow, $m_\ell$ and $o_p$.

        \item \textit{Maximum consecutive workloads.} A curriculum, $\ell$, or a professor, $p$, does not exceed their maximum allowed consecutive workload, $r_\ell$ or $s_p$. 
        
       \end{enumerate}
    
    \item \textit{Lecture assignment.} These constraints ensure that regulations regarding lecture assignments are met.
    \begin{enumerate}
        \item \textit{Single assignment.} Each lecture, $i$, is assigned to a single time slot, ($j,k$).
        
        \item\textit{Same day incompatibility.} Two lectures, $i$ and $\overline{i}$, are not taught in the same day, $j$, if the regulations state that that is the case, $\overline{i} \in d_{i}$.
        
        \item \textit{Precedence.} If there is a lecture, $\overline{i}$, which must precede another lecture, ${i}$, i.e., $\overline{i}\in e_{i}$, then, $\overline{i}$, must be taught in the same day as $i$, at an earlier hour, or a previous day. 
     
    \end{enumerate}
    \item \textit{Campus changes.} These constraints ensure that regulations regarding student or professor campus changes within the same day are met.
    \begin{enumerate}
        \item { \textit{Minimum required time.} Students of a certain curriculum, $\ell$, and professors, $p$, have enough time, $u_{\mu\overline{\mu}}$, to change campus if they have two consecutive lectures, $i$ and $\overline{i}$, in different campi, $g_i\neq g_{\overline{i}}$. }

    \end{enumerate}

\end{enumerate}

\subsection{Instances}\label{instances}
\noindent {The real-world instance, denoted as {IST\_Complete}, includes 4102 lectures, 644 classes and 900 professors. The lectures take place in 2 campi. There are lectures in the first trimester, lectures in the second trimester, and lectures that must be scheduled for the whole semester. The semester lectures must be taught in the same time slot in both trimesters. Thus, both trimesters must be scheduled simultaneously. Accordingly, each class has two distinct timetables, each corresponding to a different curriculum.
Similarly, professors and room types must have distinct timetables for each trimester. Lecture duration ranges from one to four hours, with each hour corresponding to two time slots. When there are consecutive lectures in different campi, students and professors must have at least one free time slot in between. Students and teachers can have eight hours of lectures in a single day; of those, only five hours can be taught consecutively.
}
{To expand the computational experiments, 30 random subdivisions with 200, 250, and 300 classes were obtained from instance IST\_Complete. IST has two campi where lectures are taught, with 168 and 39 rooms available for teaching. However, for the subdivisions, we considered 142 and 29 rooms for the campi. Each lecture must be assigned to a specific type of room based on its capacity. 
Instances from the 2007 International Timetabling Competition were also adapted and tested. These instances do not include multiple campi or trimesters. Moreover, maximum consecutive workload, precedence, and different day lectures were added, and unavailable slots from the original problem were not regarded.
A summary of the instances is presented in Table \ref{InstanceTable}.}

\begin{table}[h]
 \caption{Instances.}
    \label{InstanceTable}

    \begin{tabular}{|c|c|c|c|c|c|c|c|}
    \hline
 
        Instance group&Instances&Lectures ($n_i$) & Classes & Curricula ($n_\ell$)&Periods&Campi&Slots\\
        \hline
       IST\_Complete&1&4102&644&1288&2&2&24\\
       IST400&30&2200-2400&200&400&2&2&24\\
       IST500&30&2500-2700&250&500&2&2&24\\
      IST600&30&2800-3000&300&600&2&2&24\\
      ITC2007 CBCT(1-20)&20&150-450&10-150&10-150&1&1&5-9\\
       
        \hline
    \end{tabular}
 
\end{table}

{Since our instances include several characteristics that are not present in the competition instances, we decided to encode our instances with a different structure than the ones from ITC 2007. Specifically to facilitate the decomposition and increment procedures and to include other parameters such as lecture duration, campi, or trimester.}

\section{Algorithmic framework}\label{af}

\noindent In this section, we explain the algorithmic framework of our meta-heuristic. In Subsection \ref{gd}, we give a general overview. In Subsection \ref{id}, we introduce the instance decomposition and increment procedure. In Subsection \ref{cis}, we present the heuristic for generating the initial timetables. In Subsection \ref{gls}, we address the guided local search component. In Subsection \ref{alns}, we address the adaptive large neighbourhood component. In Subsection \ref{vns}, we address the variable neighbourhood search component.

\subsection {Overview}\label{gd}

\begin{figure}[h]
    \centering
    \includegraphics[scale=0.483]{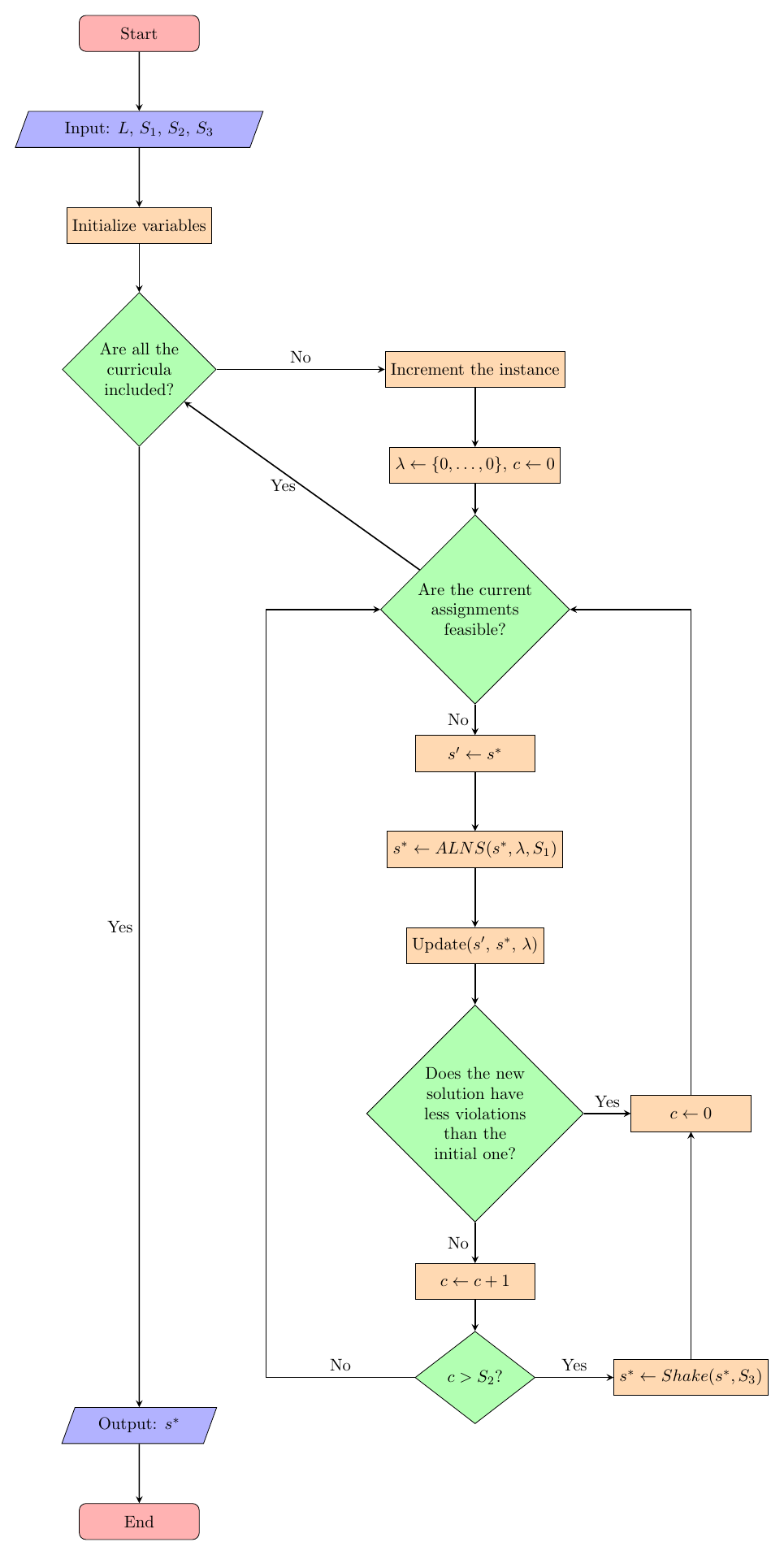}
    \caption{Meta-heuristic diagram.}
    \label{diagram overview}
\end{figure}

\noindent Our hybrid meta-heuristic uses mechanisms from adaptive large neighbourhood search (ALNS), guided local search (GLS), and variable neighbourhood search (VNS). {Its generic representation is presented in Figure \ref{diagram overview}.}
Its inputs are the full set of curricula, $L$, and three stopping criteria, $S_1$, $S_2$, and $S_3$. Additional parameters are also necessary for its functions and procedures, which are addressed in the following subsections. Its output is a feasible schedule, $s^*$.

Instead of solving the entire instance from the start, we decompose it by curricula and increment it iteratively until the set of curricula being solved, $L'$, is equal to the entire set of curricula, $L$. 
The decomposition mechanism functions in the following way: a set of lectures belonging to a set of curricula is added to the problem. Then, we find a new feasible solution, $s^*$, considering the current set of lectures, $I'$. When we have a feasible solution, we add the next set of curricula. Let us note that the assignments of each iteration can be changed if necessary in the following iterations. 

Our objective function, $f(s^*)$, represents the number of hard constraint violations in the current solution. Thus, the neighbourhoods are searched to find if there is a neighbour with fewer violations. The neighbourhood structures are explained in detail in Subsection \ref{alns}. However, instead of using a regular objective function in this search, we introduce mechanisms for penalising certain features in the solution, which are useful for escaping local optima. Accordingly, after a certain number of iterations, the feature penalty vector, $\lambda$, is updated. The features and the update procedure are explained in detail in Subsection \ref{gls}. However, sometimes this penalty is not enough to escape some local optima. To avoid those, we include a shaking phase, which occurs after a certain number of penalty recalculations without improvements. The number of local searches before updating the penalty vector is the stopping criterion $S_1$. The number of penalty weight recalculations which is attempted before inducing a shaking phase is the stopping criterion $S_2$. The duration of the shaking phase is the stopping criterion $S_3$. The pseudo-code for this meta-heuristic is presented in Algorithm \ref{alg_overview}.

\begin{algorithm}[h]
\caption{Meta-heuristic}\label{alg_overview}
	\begin{algorithmic}[1]\normalsize
	\State{\textbf{input:} $L$, $S_1$, $S_2$, $S_3$};
	\State{\textbf{output:} $s^*$};
         \State{$I' \leftarrow \{\}$};
        \State{$L' \leftarrow \{\}$};
		\State{$s^* \leftarrow \{\}$};
		\While{($length(L) > length(L')$)}
        \State{$CurriculumIncrement(s^*, L, L',I')$;}
        
        \State{$\lambda\leftarrow\{0,\ldots,0\}$;}
        \State{$c\leftarrow0$;}
        \While{($f(s^*)>0$)}
        \State{$s'\leftarrow s^*$;}
        \State{$s^* \leftarrow ALNS(s^*,\lambda,S_1)$;}
        \State{$Update(s',s^*,\lambda)$;}
        
        \If{($f(s^*) \geqslant f(s')$)}{\;$c\leftarrow c+1$;}
        \Else{\;$c\leftarrow0$;}
        \EndIf
        \If{($c> S_2$)}{\;$s^* \leftarrow Shake(s^*,S_3)$;}
        \EndIf
	    
        \EndWhile
        \EndWhile
        \State{{\textbf{return}($\mbox{$s^*$}$);}}
        \end{algorithmic} 
\end{algorithm}

\clearpage

\subsection{Instance decomposition and increment}\label{id}

\noindent {The procedure to increment the instance is straightforward. Curricula are added to the instance until a specific stopping criterion is met, or the entire set of curricula is included.}

We explore two different stopping criteria. One receives as input a fixed number of curricula to be added in each iteration, $\rho$. The other adds curricula until a certain number of hard constraint violations, $\overline{f}$, is reached. {Their respective pseudo-code algorithms are given in the Appendix. Algorithm \ref{alg_fixed_inc} presents the fixed increments procedure and Algorithm \ref{alg_vio_inc} the violations-based increment procedure. }

Whenever a new curriculum, $\ell''$, is added, its lectures, $i\in a_{\ell''}$, must be scheduled and added to the current solution, $s^*$. The function which returns the schedule for each new curriculum is presented in the following subsection.

\begin{figure}[h]
    \centering
    \includegraphics[scale=0.7]{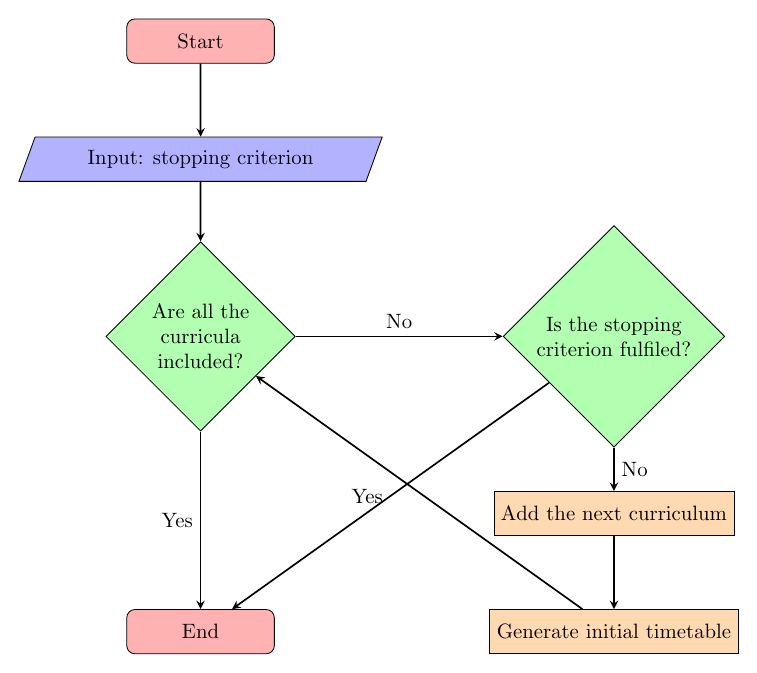}
    \caption{Curriculum increment procedure diagram.}
    \label{d-inc}
\end{figure}

Figure \ref{d-inc} provides a diagram of the procedure.
 
We test two different sorting algorithms for deciding the order in which the curricula are added. The random sorting procedure chooses a random curriculum. The ordering procedure picks a degree, then, by ascending order of year, picks a class, and then adds the first and second trimesters, each corresponding to a different curriculum, $\ell$, or timetable. This procedure is presented in Algorithm \ref{ordered_sorting}.

\begin{algorithm}[h]
\caption{Ordered sorting}\label{ordered_sorting}
	\begin{algorithmic}[1]\normalsize
	\State{\textbf{input:} $Degrees$};
	   \State{$count\leftarrow1$;}
        \For{($degree\in Degrees$)}
        \For{($year \in degree$)}
        \For{($class\in year$)}
        \For{($\ell \in class$)}
        \State{$order_\ell\leftarrow count$;}
        \State{$count\leftarrow count+1$;}
        \EndFor
        \EndFor
        \EndFor
        \EndFor
        \end{algorithmic} 
\end{algorithm}

\subsection{Initial solutions}\label{cis}

\noindent The procedure which defines the initial timetable, $ct$, for a curriculum, $\ell$, also receives as input the latest day at which the lectures can start being scheduled, $\overline{j}$, the latest hour in each day at which the lectures can start being scheduled, $\overline{k}$, the current workload for the curriculum in each day, $w'_j$, and the consecutive workload which would happen in each time-slot, ($j,k$), if a lecture, $i$, would be scheduled in such a time-slot, $r'_{ijk}$. These two parameters can be larger than 0 if the curriculum shares lectures with other already added curricula. The set of already added lectures, $I'$, is also necessary. {The procedure is summarised in Figure \ref{d-it}. The pseudo-code algorithm is presented in the Appendix in Algorithm \ref{alg_initial}. }

\begin{figure}[h]
    \centering
    \includegraphics[scale=0.7]{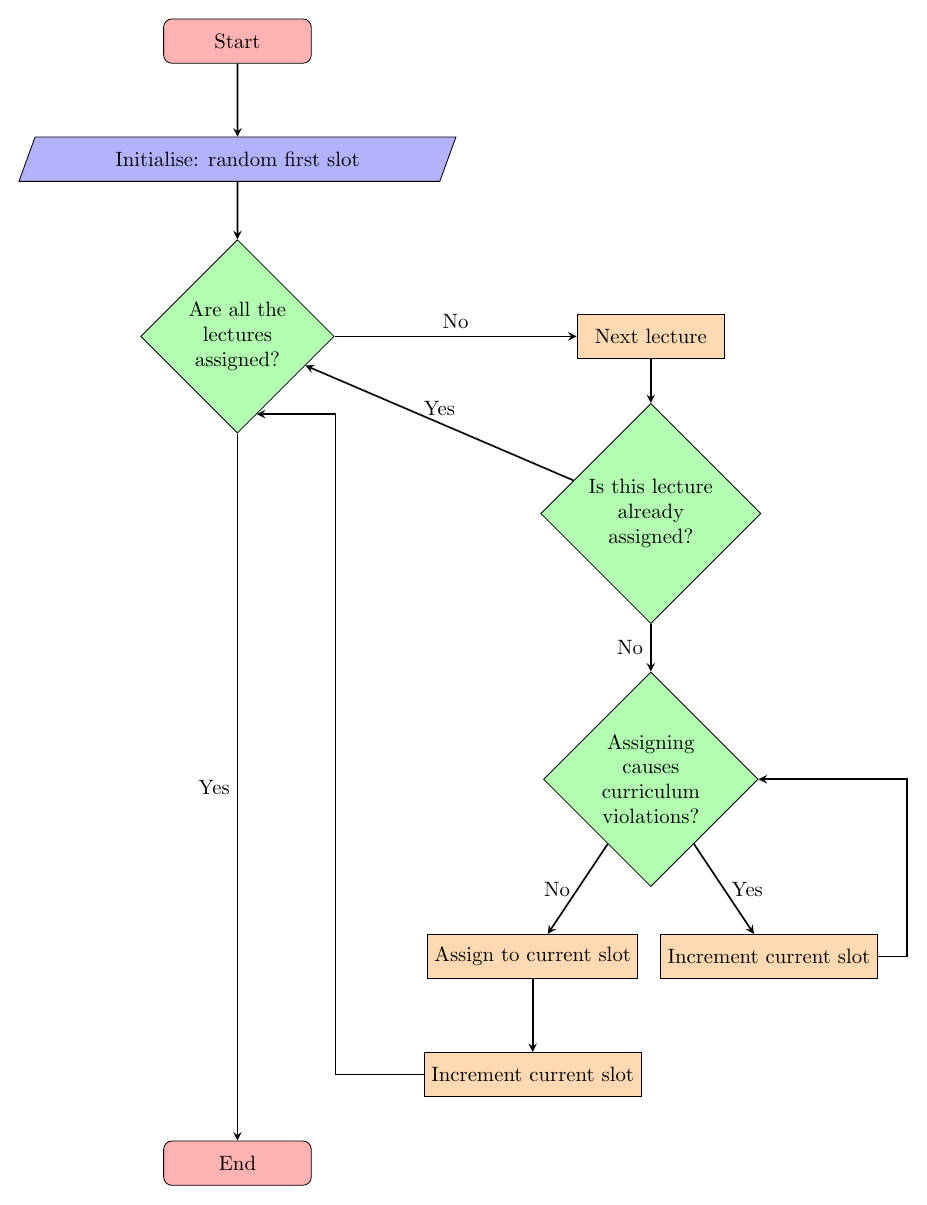}
    \caption{Initial timetables heuristic diagram.}
    \label{d-it}
\end{figure} 

This heuristic constructs the timetables to avoid maximum daily and consecutive workload and juxtaposition violations for curricula. Nonetheless, some of these violations may be present in the initial timetable if lectures it shares with other curricula had already been scheduled in a way that violates those constraints. 

For each lecture in the curriculum that has not yet been added, $i\in a_{\ell} \land i \notin I'$, we check if scheduling it for the current time slot would cause any of the aforementioned violations. If it would, we go to the next time slot. When a suitable time slot is found, the lecture is assigned, and the daily, $w'_j$, and consecutive, $r'_{ijk}$, workload parameters are updated accordingly. 

\clearpage

\subsection{Guided local search}\label{gls}

\noindent In GLS, the objective function, $f(s^*)$, is augmented based on certain solution features, $\psi(s^*)$. This can be useful to escape local optima by penalising the features present in such solutions. The features we implemented are of three types: the constraint violations for specific curricula, $\psi^1_\ell$; the total violations for all professors, $\psi^2$, all rooms, $\psi^3$, all precedence rules, $\psi^4$, and all different day rules, $\psi^5$; and the violations occurring in specific time slots, $\psi^6_{jk}$. However, these features are not always penalised. A vector, $\lambda$, takes the value 0 for non-penalised features and a value greater than 0 otherwise. A larger value for a certain feature means it is more penalised in the objective function. 

Considering the notation: $ f_{jk\ell}$, are the constraint violations for a curriculum, $\ell$, in a time slot, ($j,k$); $ f_{jkp}$, are the constraint violations for a professor, $p$, in a time slot, ($j,k$); $ f_{jkt\mu}$, are the constraint violations for a room type, $t$, in a campus, $\mu$, in a time slot, ($j,k$); $ f^4_{ijk}$, are the precedence constraint violations for a lecture, $i$, in a time slot, ($j,k$); and $ f^5_{ijk}$, are the different day constraint violations for a lecture, $i$, in a time slot, ($j,k$). Objective function \eqref{objective_function} functions similarly to an objective function in any minimisation problem. Nonetheless, in our problem, when it reaches the value of 0, it signifies that a feasible solution is found. 

\begin{multline}\label{objective_function}
    \displaystyle \min f(s^*) = \sum_{j=1}^{n_j}\sum_{k=1}^{n_k}\sum_{\ell=1}^{n_\ell}\left( f_{jk\ell} + \lambda^1_\ell\psi^1_{jk\ell} + \lambda^6_{jk}\psi^6_{jk} (f_{jk\ell}>0)\right) + \\
    \sum_{j=1}^{n_j}\sum_{k=1}^{n_k}\sum_{p=1}^{n_p}\left( f_{jkp} + \lambda^2\psi^2_{jkp}+\lambda^6_{jk}\psi^6_{jk} (f_{jk\ell}>0)\right) + \\
    \sum_{j=1}^{n_j}\sum_{k=1}^{n_k}\sum_{t=1}^{n_t}\sum_{\mu=1}^{n_\mu}\left( f_{jkt\mu}+\lambda^3\psi^3_{jkt\mu} + \lambda^6_{jk}\psi^6_{jk} (f_{jkt\mu}>0)\right) + \\
    \sum_{i=1}^{n_i}\sum_{j=1}^{n_j}\sum_{k=1}^{n_k}\left( f^4_{ijk}+ \lambda^4\psi^4_{ijk} + f^5_{ijk}+ \lambda^5\psi^5_{ijk} + \lambda^6_{jk}\psi^6_{jk}\left( (f^4_{ijk}>0) + (f^5_{ijk}>0)\right)\right)
\end{multline}

The general algorithm framework is presented in Algorithm \ref{gls-a}.

\begin{algorithm}[h]
\caption{Guided local search general framework: \textit{GLS()}}\label{gls-a}
\begin{algorithmic}[1]\normalsize
	\State{\textbf{input:} $s^*$, $\lambda$};
	\State{\textbf{output:} $s^*$};
        
        \While{(stopping criterion $=false$)}
       
        \State{$s'\leftarrow Search(s^*)$};
        \If{($f(s') + \lambda\psi(s')\leqslant f(s^*)+ \lambda\psi(s^*)$)}{\;$s^*\leftarrow s'$};
        \EndIf
        \State{$Update(\lambda)$};
        \EndWhile
        \State{{\textbf{return}($\mbox{$s^*$}$);}}
        \end{algorithmic} 
\end{algorithm}

The procedure to update the values of the penalty vector, $\lambda$, also considers an initial solution, $s'$, the current solution, $s^*$, and the set of features, $\beta$. 

If the current solution, $s^*$, performs better than the initial solution, $s'$, the penalty vector, $\lambda$, is set to 0 for every feature. Otherwise, we assume we are in a local optimum for the current penalty vector, $\lambda$. Accordingly, for every feature, $i\in\beta$, which has an equal or worse performance in the current solution, $\psi^i(s^*)$, the penalty, $\lambda_i$, is incremented by 1. Otherwise, for a feature that performs better in the current solution or is not present, the penalty is set to 0. 

Specifically, for the curricula feature type, if a curriculum has the same or more violations in the current solution, the penalty for such violations is incremented. If the value of the total professor violations is greater or equal in the current solution, the penalty for such violations is incremented, and the same for the room, precedence, and different day violations. If there are the same or more violations in a specific time slot in the current solution, the penalty for any violation occurring in that time slot is incremented. This procedure is presented in Algorithm \ref{update_gls}.

\begin{algorithm}[h]
\caption{ \textit{Update()}}\label{update_gls}

	\begin{algorithmic}[1]\normalsize

	\State{\textbf{input:} $s'$, $s^*$, $\lambda$, $\beta$};

        \If{($f(s^*)+\lambda\psi(s^*)\geqslant f(s')+\lambda\psi(s')$)}
        \For{($i \in \beta$)}
             \If{($0<\psi_i(s^*)\geqslant \psi_i(s')$)}
                \State{$\lambda_i\leftarrow \lambda_i+1$};
                \Else
                \State{$\lambda_i\leftarrow 0$};
             \EndIf
             
        \EndFor
        \Else
        \State{$\lambda=\{0,\ldots,0\}$};
        \EndIf
        \end{algorithmic} 
\end{algorithm}

Considering each type of constraint violation as an objective, this methodology can be analysed as a multi-objective search technique. Effectively, what is happening is that out of the objective functions being regarded, the most challenging to improve are identified. Then, their weights, $\lambda_i$, are updated to help the search move towards non-dominated solutions that improve these objectives' performance, even if the total sum of the effective violations increases.

\subsection{Adaptive large neighbourhood search}\label{alns}
\noindent In ALNS, the local search is based on different neighbourhood structures. In each iteration, a neighbourhood structure, $n$, is selected based on a probability vector, $p$, and a local search is conducted. If a better solution is found, it replaces the best-known solution. Let us note that we use an augmented objective function adapted from GLS procedures, $f'(s)$. Then, based on the current solution characteristics and the previous successes of searches in each neighbourhood, the probability vector, $p$, is updated. {The general algorithm framework is presented in Algorithm \ref{alg_alns}.}

\begin{algorithm}[h]
\caption{Adaptive large neighbourhood search general framework: \textit{ALNS()}}\label{alg_alns}
	\begin{algorithmic}[1]\normalsize
	\State{\textbf{input:} $s^*$,$\lambda$, $N$, $p$, $S_1$};
	\State{\textbf{output:} $s^*$};
        
        \For{($iteration=1,\ldots,S_1$)}
        \State{$n\leftarrow SelectNeighbourhood(N,p)$};
        \State{$s'\leftarrow Search(s^*,n)$};
        \If{($f'(s')\leqslant f'(s^*)$)}{\;$s^*\leftarrow s'$};
        \EndIf
        \State{$Update(p)$};
        \EndFor
        \State{{\textbf{return}($\mbox{$s^*$}$);}}
        \end{algorithmic} 
	   
\end{algorithm}

To implement ALNS, it is necessary to define the set of neighbourhood structures, $N$, the procedure for updating the probability vector, $p$, and the stopping criterion. We established six main neighbourhood structures, which we named: worst slot, worst curriculum, worst professor, worst room type, different day, and precedence neighbourhoods. Moreover, each of them has two sub-types, named swap and non-swap. In each neighbourhood, a lecture is selected, and then, in the swap sub-type, its time slot is swapped with another lecture, and in the non-swap sub-type, it is randomly assigned to another time slot.

\begin{figure}[h]
    \centering
    \includegraphics[scale=0.8]{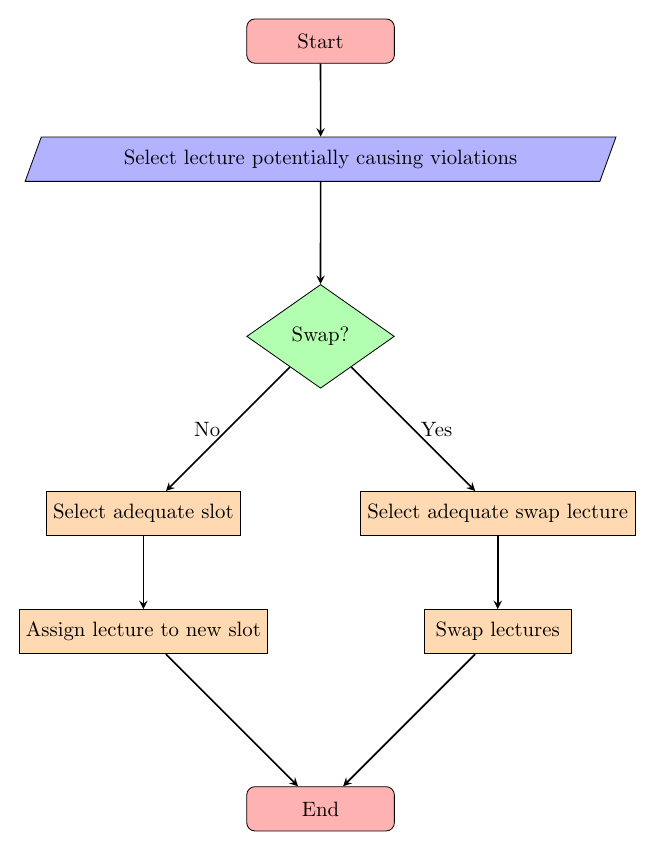}
    \caption{Neighbourhood structures general framework.}
    \label{d-ns}
\end{figure}

{The neighbourhoods follow a general framework that can be applied to each type of constraint violation. The first step is randomly selecting a lecture potentially causing violations. Then, either an adequate new time slot or an adequate swap lecture is randomly assigned. This adequacy is assessed based on the type of violation being caused. This general framework is presented in Figure \ref{d-ns}.}

\clearpage

Specifically, each neighbourhood structure operates in the following manner:

\begin{itemize}[label={--}]
    \item \textit{Worst slot}. A lecture, $i'$, is selected from the set of lectures assigned during the time slot with the most violations, $I''$. In the non-swap sub-type, it is assigned to a random time slot. In the swap sub-type, a random curriculum, $\ell'$, is selected from the set of curricula of the searched lecture, $i'$. Then, a different lecture, $i''$, is selected from this curriculum's set of lectures, $a_{\ell'}$, for a swap. {This neighbourhood structure is presented in the Appendix in Algorithm \ref{alg_worst_slot}.}
    \item \textit{Worst curriculum}. A curriculum, $\ell'$, is selected from the set of the curricula with the most violations, $L''$. Then, a lecture belonging to that curriculum, $i' \in a_{\ell'}$, is selected from the set of lectures assigned during the time slot with the most violations for that curriculum, $I''$. In the non-swap sub-type, it is assigned to a random time slot. In the swap sub-type, a random lecture, $i''$, is selected from this curriculum's set of lectures, $a_{\ell'}$, for a swap. {This neighbourhood structure is presented in the Appendix in Algorithm \ref{alg_worst_curr}.}
    \item \textit{Worst professor}. A professor, $p'$, is selected from the set of professors with the most violations, $P''$. Then, a lecture taught by that professor, $i' \in a_{\ell'}$, is selected from the set of lectures assigned during the time slot with the most violations for the professor, $I''$. In the non-swap sub-type, it is assigned to a random time slot. In the swap sub-type, a random curriculum, $\ell'$, is selected from the set of curricula of the searched lecture, $i'$. Then, a different lecture, $i''$, is selected from this curriculum's set of lectures, $a_{\ell'}$, for a swap.{ This neighbourhood structure is presented in the Appendix in Algorithm \ref{alg_worst_prof}.}
    \item \textit{Worst room type}. A lecture, $i'$, taught in the room type with the most violations, $v_{i't}=1$, in the campus with the most violations, $\mu$, during the time slot with the most violations for that campus and room type, ($\hat{j}, \hat{k}$) is selected from the set of lectures in a similar situation, $I''$. In the non-swap sub-type, it is assigned to a random time slot. In the swap sub-type, a random curriculum, $\ell'$, is selected from the set of curricula of the searched lecture, $i'$. Then, a different lecture, $i''$, is selected from this curriculum's set of lectures, $a_{\ell'}$, for a swap.{ This neighbourhood structure is presented in the Appendix in Algorithm \ref{alg_worst_room}.}
    \item \textit{Different day}. A lecture, $i'$, that violates the different day constraints is selected. In the non-swap sub-type, it is assigned to a random time slot, on a different day. In the swap sub-type, a random curriculum, $\ell'$, is selected from the set of curricula of the searched lecture, $i'$. Then, a different lecture, $i''$, assigned to a different day, is selected from this curriculum's set of lectures, $a_{\ell'}$, for a swap. {This neighbourhood structure is presented in the Appendix in Algorithm \ref{alg_diff_da}.}
    \item \textit{Precedence}. A lecture, $i'$, and a different lecture, $i''$, which should be scheduled after the first and is not, are selected. In the non-swap sub-type, the first lecture, $i'$, is assigned to a random time slot, before the time slot in which the other lecture is scheduled. In the swap sub-type, both lectures are swapped. {This neighbourhood structure is presented in Appendix in Algorithm \ref{alg_prec}.}
\end{itemize}

Except for the worst slot neighbourhood, each other neighbourhood structure, $n$, targets a specific type of violation, $f_n(s^*)$. Accordingly, in each ALNS iteration, the probability of choosing a neighbourhood structure, $p_n$, is updated based on the number of violations of each type, $f_n(s^*)$, and the number of times a search in a neighbourhood was successful at finding an improved solution, $\Delta_n$. The update procedure is presented in Algorithm \ref{update_alns}.

\begin{algorithm}[h]
\caption{\textit{Update()}}\label{update_alns}
\begin{algorithmic}[1]\normalsize
	\State{\textbf{input:} $N$, $\Delta$, $s^*$};

        \State{$p'_1\leftarrow \Delta_1 avg(f(s^*_n))$};
        
        \For{($n \in \{2, \ldots, 6\}$)}
        \State{$p'_n\leftarrow \Delta_n f(s^*_n)$};
        \EndFor

         \For{($n \in N$)}
        \State{$p_n\leftarrow \frac{p'_n}{\sum_{n=1}^6p'n}$};
        \EndFor
        
        \end{algorithmic} 
	   
\end{algorithm}

The stopping criterion for our ALNS implementation is based on a maximum number of iterations. This means that the number of times a search in a neighbourhood successfully found an improved solution, $\Delta_n$, is reset after this number of iterations.

\subsection{Variable neighbourhood search}\label{vns}

\noindent In VNS, whenever local optima are identified, a shaking phase is used to move to different areas of the solution space. In our meta-heuristic, such a shaking phase is conducted if no improvement moves are found after a certain number of feature penalty updates are made. This shaking phase functions similarly to the ALNS algorithm, but every move is accepted instead of only those which introduce improvements. This phase is stopped after a number of iterations, $S_3$. It is presented in Algorithm \ref{shakee}.

\begin{algorithm}[h]
\caption{\textit{Shake()}}\label{shakee}
	\begin{algorithmic}[1]\normalsize
	\State{\textbf{input:} $s^*$, $N$, $p$, $S_3$};
	\State{\textbf{output:} $s^*$};
        
        \For{($iteration=1,\ldots,S_3$)}
        \State{$n\leftarrow SelectNeighbourhood(N,p)$};
        \State{$s^*\leftarrow Search(s^*,n)$};
        
        \State{$Update(p)$};
        \EndFor
        \State{{\textbf{return}($\mbox{$s^*$}$);}}
        \end{algorithmic} 
	   
\end{algorithm}

\section{Computational experiments}\label{ce}

\noindent In this section, we address the computational experiments. In Subsection \ref{Scp}, we test several stopping criteria parameters for smaller instances. In Subsection \ref{IDF}, we present the results of instance decomposition with fixed curriculum increments for smaller instances. In Subsection \ref{IDV}, we present the results of instance decomposition with violations based curriculum increments for smaller instances. In Subsection \ref{EL}, we test the meta-heuristic with larger subdivisions and the real-world instance. In Subsection \ref{Benchmark}, we present the experiments with the adapted ITC 2007 CBCT instances.

All algorithms are implemented in C++. For the IST instances, the CPU is an 11th Gen Intel(R) Core(TM) i5-11400 @ 2.60GHz   2.59 GHz, and the installed RAM is 16.0 GB (15.8 GB usable). For the adapted competition instances, the CPU is an Intel(R) Core(TM) i7-8565U CPU @ 1.80GHz   1.99 GHz, and the installed RAM is 8.00 GB (7.78 GB usable).

\subsection{Stopping criteria parameters}\label{Scp}

\noindent In this subsection, we test the effect of three parameters when solving instances IST400 without instance decomposition. Specifically:

\begin{itemize}[label={--}]

\item $S_1$. Number of ALNS iterations without improvement before recalculating the GLS penalty weights;

\item $S_2$. Number of penalty recalculations without improvement before triggering a shaking phase;

\item $S_3$. Duration of the shaking phase.
\end{itemize}

\begin{figure}[h]
    \centering
    \includegraphics[scale=0.2]{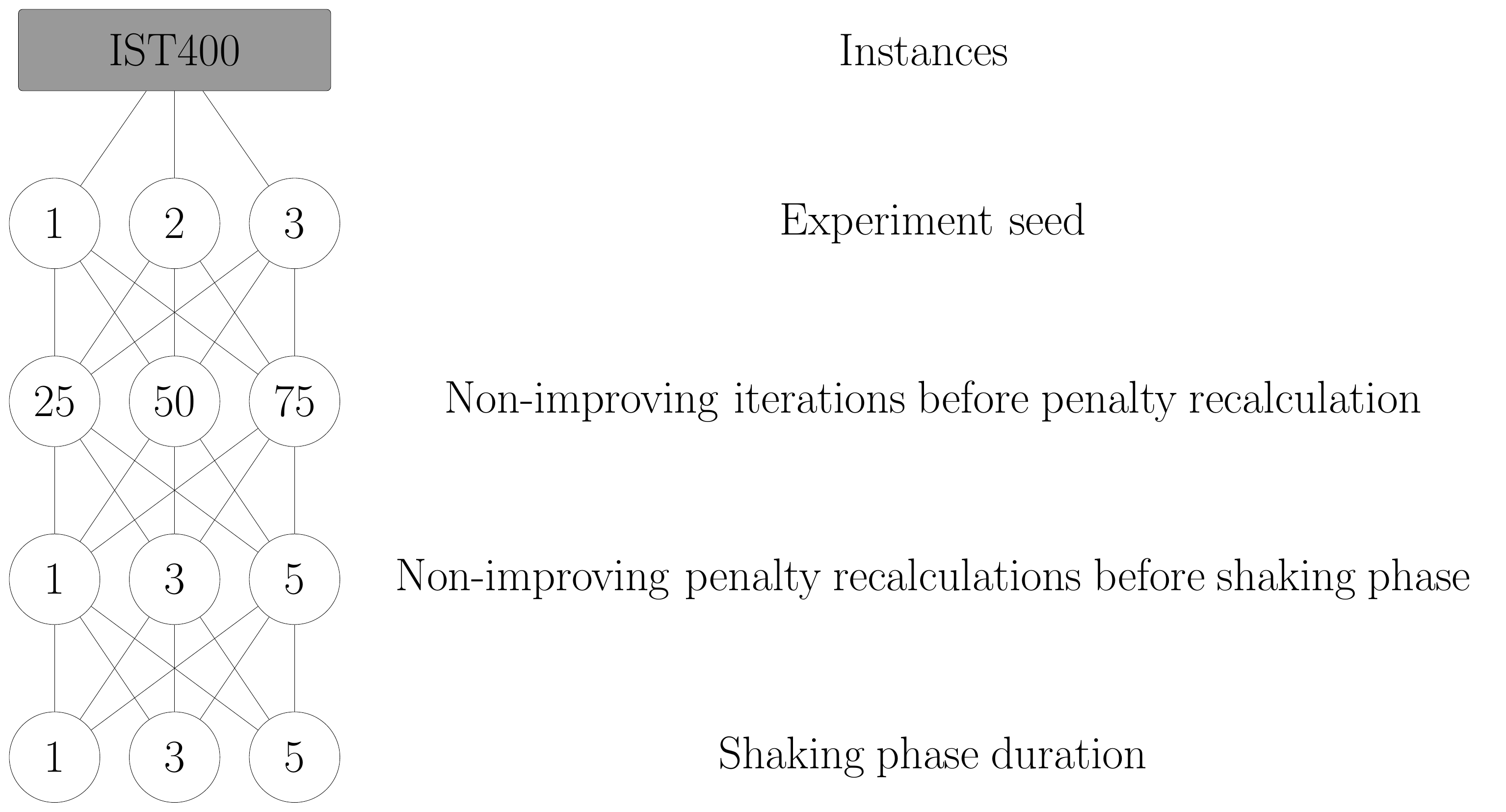}
    \caption{Stopping criteria parameter experiments ($S_1,S_2,S_3$).}
    \label{exp-Scp}
\end{figure}

Parameter combinations with a larger $S_1$ allow for greater intensification of a specific search area before inducing the penalty weight recalculation, a diversification mechanism. Changing these penalties changes the trade-off profile between different infeasible solutions. For example, if there has been a time slot that has had a large number of constraint violations for the last iterations, and that number has not improved, the weight for violations on that time slot is increased, and perhaps a new solution can be found with fewer violations on that time slot, even if it has more on another. A shaking phase is introduced if these mechanisms fail to diversify the search to escape a local optimum for $S_2$ recalculations. During this phase, $S_3$ iterations are made without checking if the new solution is better than the previous one. Accordingly, parameter combinations with smaller $S_1$ and $S_2$ are more active in diversifying the search but less precise in identifying local optima. Smaller shaking phases, $S_3$, introduce fewer violations in the new solutions but might not introduce enough change to escape local optima. Each instance and parameter combination is tested with three different seeds to reduce randomness effects. The experiments are summarised in Figure \ref{exp-Scp}.

\clearpage

The average solution time across all 30 instances and for 3 runs of each instance and configuration is presented in Table \ref{sotime-scp}. The five fastest combinations are in bold.

\begin{table}[h]
  \caption{Solution time for each stopping criteria combination.}
  \label{sotime-scp}
  \centering
  \begin{tabular}{@{}ccc@{}}
    \parbox{0.3\textwidth}{
      \centering
    
       \begin{tabular}{|c|c|c|c|}
    \hline
     $S_1$&$S_2$&$S_3$&CPU(s)  \\
     \hline
      25 & 1&1&140\\
      25&1&3&157\\
      25&1&5&206\\

      \textbf{25}&\textbf{3}&\textbf{1}&	\textbf{125}\\
    \textbf{25}&\textbf{3}&\textbf{3}	&\textbf{124}\\
\textbf{25}&\textbf{3}&\textbf{5}&	\textbf{126}\\
25&5&1&	128\\
25&5&3	&131\\
25&5&5&	130\\
\hline
\end{tabular}
    }
    &
    \parbox{0.3\textwidth}{
      \centering
     
       \begin{tabular}{|c|c|c|c|}
 \hline
     $S_1$&$S_2$&$S_3$&CPU(s)  \\
     \hline
50&1&1&	129\\
\textbf{50}&\textbf{1}&\textbf{3}&	\textbf{125}\\
50&1&5&	137\\
50&3&1&	133\\
50&3&3&	139\\
50&3&5&	134\\
50&5&1&	138\\
50&5&3&	136\\
50&5&5&	138\\
\hline
\end{tabular}
    }
    &
    \parbox{0.3\textwidth}{
      \centering
     
      \begin{tabular}{|c|c|c|c|}
 \hline
     $S_1$&$S_2$&$S_3$&CPU(s)  \\
     \hline
75&1&1&	133\\
75&1&3&	129\\
75&1&5&	134\\
75&3&1&	134\\
75&3&3&	127\\
\textbf{75}&\textbf{3}&\textbf{5}&	\textbf{126}\\
75&5&1&	139\\
75&5&3&	130\\
75&5&5&	130\\

      \hline
    \end{tabular}
    }
  \end{tabular}
\end{table}

\begin{figure}[h]
    \centering
\begin{tikzpicture}[scale=0.7]
	\begin{axis}[%
 xlabel={($S_1,S_2$)},
	ylabel={CPU(s)}
	={(xy): \thisrow{label}},%
 	xtick = {1,2,3,4,5,6},
		xticklabel style = {align=center, font=\normalsize},
		xticklabels = {(25,1), (50,1), (75,1), (25,5), (50,5),  (75,5)},
	scatter/classes={%
		a={mark=square*,black},%
		b={mark=square*,black}}    ]
	\addplot[scatter,only marks,%
		scatter src=explicit symbolic]%
	table[meta=label] {
x  y   label
1  140 a
1  157 a
1  206 a

2  129 a
2  125 a
2  137 a

3  133 a
3  129 a
3  134 a

4  128 b
4  131 b
4  130 b

5  138 b
5  136 b
5  138 b

6  139 b
6  130 b
6  130 b
	};
 
	\end{axis}
\end{tikzpicture}
 \caption{Active \textit{vs} precise parameter combinations.}
    \label{DEG performance}
\end{figure}
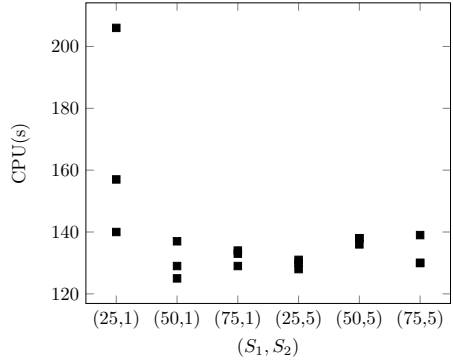

Considering the activity \textit{vs} precision framework, more active parameter combinations' average solution times are more affected by the values of other parameters as they perform more diversification moves. For $S_1=25$, we have the fastest (25,3,3)=124 seconds and the slowest combination (25,1,5)=206s, which have an 82 seconds solution time difference. Contrarily, the difference between the fastest and slowest combinations for $S_1=50$ is 14 seconds, and $S_1=75$ is 12 seconds. The same conclusion can be drawn for $S_2$ and combinations of $S_1$ and $S_2$. Figure \ref{DEG performance} represents this effect, showing a gradual decrease in average time variation. However, the (75,5) parameter combination breaks this trend. (75,5,1) is 9 seconds slower than the other (75,5) combinations, which have the same average solution time. This might be a statistical outlier with no significant meaning, or it might indicate that smaller shaking phase durations, $S_3$, for increasingly precise combinations can make certain local optima considerably harder to overcome.

This effect does not mean that more active combinations do not perform very well. Four out of the five fastest combinations can be considered relatively active: (25,3,1), (25,3,3), (25,3,5), and (50,1,3).

\subsection{Instance decomposition with fixed curriculum increments}\label{IDF}

\noindent In the previous subsection, the meta-heuristic received the data for all the curricula from the start. In this subsection, we address the experiments regarding fixed curricula increments. The mechanism tested here works as follows: a set of lectures belonging to a fixed number of curricula, $\rho$, is added to the problem. Then, we find a new feasible solution, $s^*$, considering the current set of lectures, $I'$. When we have a feasible solution, we add the next set of $\rho$ curricula until the complete set is solved. The assignments of each iteration can be changed if necessary in the following iterations. 

These experiments test combinations of three factors in the solution times of instances IST400:

\begin{itemize}[label={--}]

\item \textit{Curriculum order}. Randomised curricula increments and ordered increments, following the order from Algorithm 4, which clusters classes based on degrees and years. This means that curricula with more lectures and teachers in common are grouped;

\item \textit{Parameterisation}. Impact of more active (25,3,1) \textit{vs} more precise (75,3,5) parameter combinations;

\item \textit{Curricula added per iteration}. Number of curricula added in each iteration.

\end{itemize}

\begin{figure}[h]
    \centering
    \includegraphics[scale=0.2]{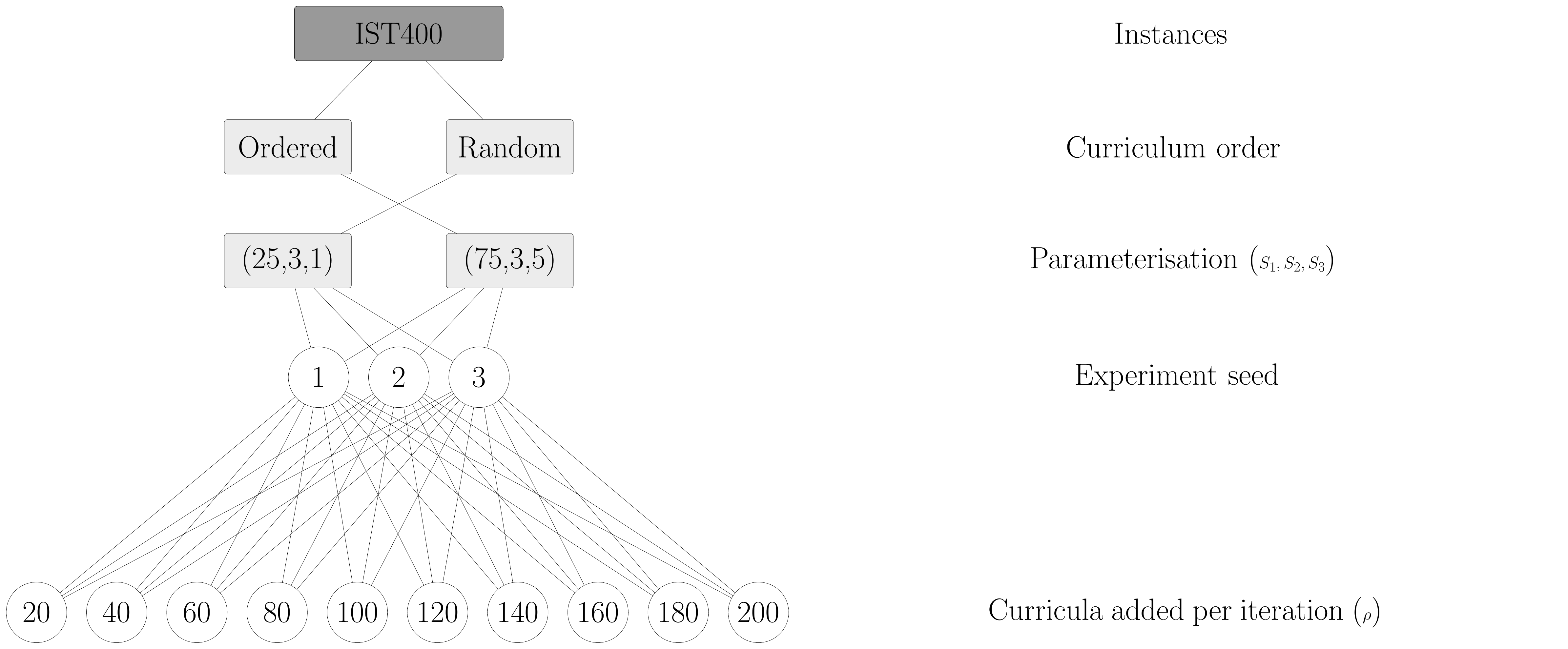}
    \caption{Fixed curricula increment experiments.}
    \label{exp_fixed_curr_400}
\end{figure}

The experiments are summarised in Figure \ref{exp_fixed_curr_400}. 

\begin{figure}[h]
    \centering
\begin{tikzpicture}[scale=0.7]
	\begin{axis}[%
  xmin=0,xmax=220,ymin=85,ymax=160,
  xtick distance = 40, 
  ytick distance = 10, 
 xlabel={Curricula added per iteration},
	ylabel={CPU(s)}]
\addplot coordinates {
	(020,	106)
(040,	97)
(060,	94)
(080,	91)
(0100,	97)
(0120,	94)
(0140,	95)
(0160,	94)
(0180,	98)
(0200,	96)

};

\addplot coordinates {
(20,	119)
(40	,109)
(60	,111)
(80	,108)
(100,	102)
(120,	108)
(140,	104)
(160,	106)
(180,	108)
(200,	103)

};

\addplot coordinates {
(020,	126)
(040,	117)
(060,	119)
(080,	115)
(0100,	114)
(0120,	117)
(0140,	113)
(0160,	117)
(0180,	117)
(0200,	117)

};

\addplot+[no marks] coordinates {
(020,	124)
(040,	124)
(060,	124)
(080,	124)
(0100,	124)
(0120,	124)
(0140,	124)
(0160,	124)
(0180,	124)
(0200,	124)

};

\legend{{Ordered, (25,3,1)},{Ordered, (75,3,5)},{Random, (25,3,1)},Best without decomposition }
	\end{axis}
\end{tikzpicture}
 \caption{Fixed curricula increments experiment results.}
    \label{exp_results_fd}
\end{figure}
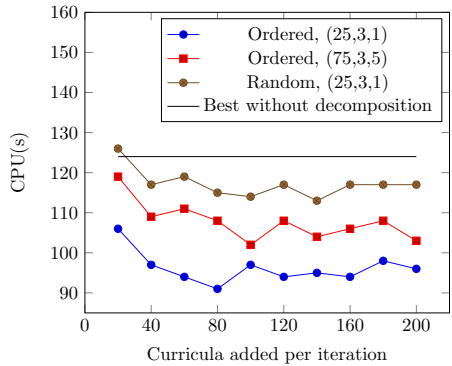

The results show that decomposing and incrementing the problem can reduce solution times to a significant degree. Considering the best curricula added per iteration value, $\rho$, for each group of experiments: the experiments ``Ordered, (25,3,1), $\rho=80$'' revealed a solution time reduction of 27.2\% when compared to the same parameterisation, but without decomposition; the experiments ``Ordered, (75,3,5), $\rho=100$'', promoted a 20.6\% solution time reduction; and the experiments ``Random, (25,3,1), $\rho=140$'', promoted a 9.6\% solution time reduction. These results are presented in Figure \ref{exp_results_fd}.

For smaller $\rho$ values, a decrease in this parameter indicates a trend for an exponential increase in solution times. This is an expected behaviour, as if the number of curricula added per iteration is too small, the probability that their timetables might have to be changed in future iterations increases. 

The order in which the curricula are added also greatly impacts the efficiency of the methodology, as evidenced by both ordered experiments outperforming the random experiments for every $\rho$ value. The ordering mechanism clusters curricula based on their degrees and years. These results indicate that solving curricula with more classes and professors in common at the same time is beneficial.

To reduce the amount of time the meta-heuristic takes to adjust schedules that will likely need further adjustment in the future, the increments need to be larger. This is especially important when the curricula are not grouped together in the random experiments. Accordingly, these performed better at larger $\rho$ values, 140. Contrarily, the ordered experiments performed better in the 80-100 range. In the ordered experiments, the curricula added in each iteration are less likely to introduce violations based on the assignments of already solved curricula. For example, let us say that three classes from the same year and degree (A, B, C) have to choose two out of three elective courses (1, 2, 3). Class A are the students who choose courses 1 and 2, class B choose 1 and 3, and class C choose 2 and 3. In a random increment, classes A and B might be solved initially, and their entire timetables are built around having lectures from courses 2 and 3 at similar time slots. When class 3 is added later on, these assignments are now infeasible and might be difficult to clear up, as to remove them, some moves which introduce other violations might be necessary. Contrarily, if the increments are made in an ordered manner, this type of situation is less likely to occur.

The parameterisation which is more active in introducing diversification mechanisms (25,3,1), performed better than the parameterisation which intensifies the search more within the same area before diversifying (75,3,5). These combinations are very similar when solving the problems without decomposition, taking, on average, 125s and 126s, respectively. When using decomposition, there are always fewer violations and fewer curricula, professors, and rooms have infeasible schedules during each iteration. This means that local optima are going to be easier to identify. Thus, the extra intensification of a certain area promoted by the combination (75,3,5) is not so beneficial when dealing with a decomposed problem. In this case, active combinations are more efficient, as performing more intensification iterations is less likely to find better solutions.

\subsection{Instance decomposition with violations based curriculum increments}\label{IDV}

\noindent In this subsection, we present the results regarding violations-based curricula increments. The mechanism tested here works as follows: a set of lectures belonging to one curriculum is added to the problem. Then, if the number of constraint violations surpasses parameter $\overline{f}$, a new feasible solution, $s^*$, is found considering the current set of lectures, $I'$. This process is repeated until the complete set of curricula is solved. The assignments of each iteration can be changed if necessary in the following iterations. 

\begin{figure}[h]
    \centering
    \includegraphics[scale=0.2]{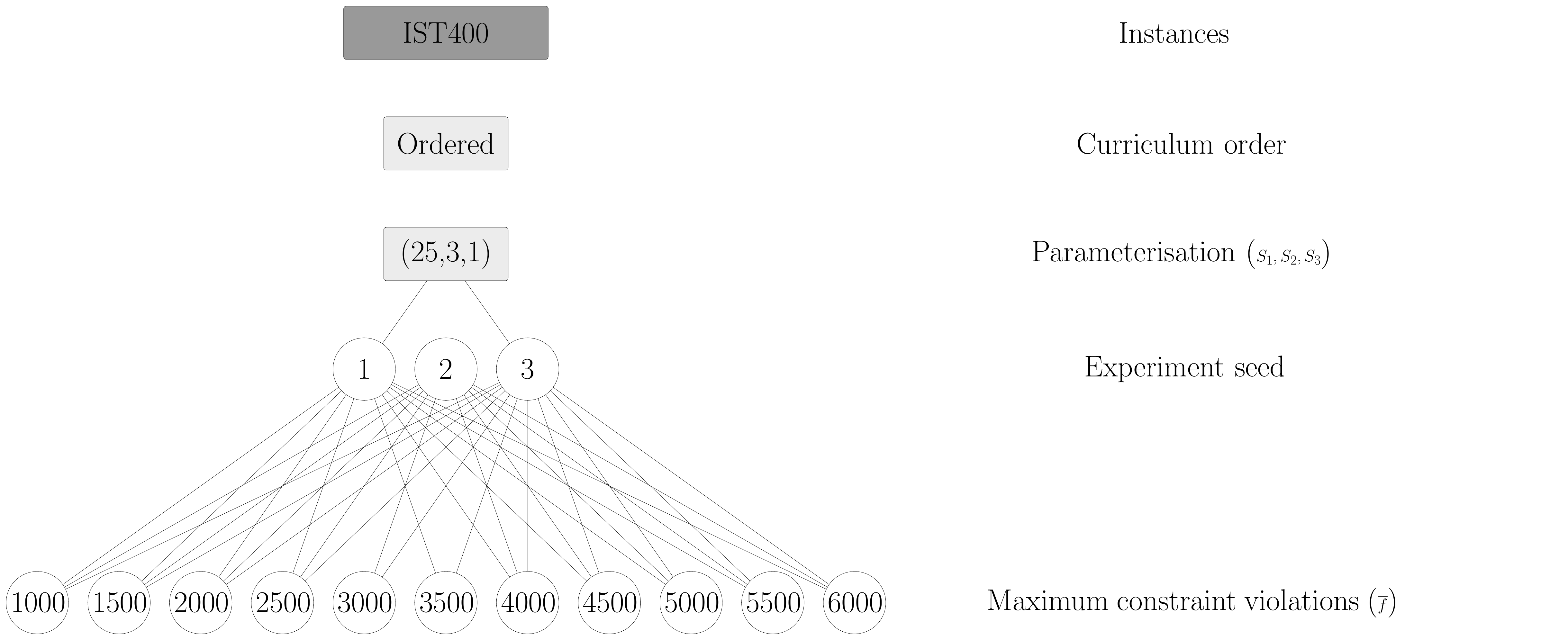}
    \caption{Violations-based curricula increment experiments.}
    \label{exp_vio}
\end{figure}

These experiments test the effect of altering the maximum number of violations before the start of the solution of the current set of lectures, $I'$, on instances IST400. They are summarised in Figure \ref{exp_vio}.

As a reference, 500 violations represent roughly 12.5 curricula. Nonetheless, this value varies quite significantly. Table \ref{equivalent curricula} provides reference values for this estimation. 

\begin{table}[h]
    \centering
     \caption{Estimate of the curricula added per iteration based on the maximum violations.}
 
    \begin{tabular}{|c|c|}
        \hline
         Violations& Estimated curricula added \\
         \hline
         1000	&	25 \\
1500	&	38\\
2000	&	50\\
2500	&	63\\
3000	&	75\\
3500	&	88\\
4000	&	100\\
4500	&	113\\
5000	&	125\\
5500	&	138\\
6000	&	150\\
 \hline
    \end{tabular}
   
    \label{equivalent curricula}
\end{table}

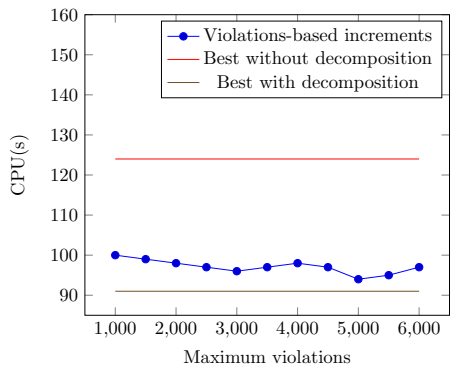
\begin{figure}[h]
    \centering
\begin{tikzpicture}[scale=0.7]
	\begin{axis}[%
  xmin=500,xmax=6500,ymin=85,ymax=160,
  xtick distance = 1000, 
  ytick distance = 10, 
 xlabel={Maximum violations},
	ylabel={CPU(s)}]
\addplot coordinates {
	(1000,	100)
(1500,	99)
(2000,	98)
(2500,	97)
(3000,	96)
(3500,	97)
(4000,	98)
(4500,	97)
(5000,	94)
(5500,	95)
(6000,	97)

};

\addplot+[no marks] coordinates {
(1000,	124)
(1500,	124)
(2000,	124)
(2500,	124)
(3000,	124)
(3500,	124)
(4000,	124)
(4500,	124)
(5000,	124)
(5500,	124)
(6000,	124)

};

\addplot+[no marks] coordinates {
(1000,	91)
(1500,	91)
(2000,	91)
(2500,	91)
(3000,	91)
(3500,	91)
(4000,	91)
(4500,	91)
(5000,	91)
(5500,	91)
(6000,	91)

};

\legend{{Violations-based increments},Best without decomposition, Best with decomposition }
	\end{axis}
\end{tikzpicture}
 \caption{Violations-based curricula increments experiment results.}
    \label{exp_results_vd}
\end{figure}

{This mechanism requires the constant generation and assessment of the assignments for each curriculum being added. Comparatively, when doing fixed increments this process is not required}. Nonetheless, it presents a trade-off by being more flexible in the number of curricula being added, and doing smaller increments when a feasible solution is expected to be more challenging to find. Figure \ref{exp_results_vd} compares these experiments and the best results of previous experiments.

In these experiments, the trade-off of the added flexibility proved insufficient to overcome the added requirements of assessing the violations during the increment phase, and a better average solution time could not be found. Still, the mechanism proved competitive with the previous increment procedure, especially for larger values of estimated curricula added. For the experiments with fewer than 100 estimated curricula added, the violations-based experiments took, on average, between 2 and 5 more seconds than their closest comparison regarding the fixed increment experiments. However, for larger maximum violations limits, they perform very similarly: $\rho=100, CPU=97s$ and $\overline{f}=4000, CPU=98s$, $\rho=120, CPU=94s$ and $\overline{f}=5000, CPU=94s$, and $\rho=140, CPU=95s$ and $\overline{f}=5500, CPU=95s$. This indicates that for larger fixed increments, where the absolute variation of the number of violations is potentially greater, the added flexibility of this mechanism starts to overcome the extra assessment requirements.

\subsection{Experiments with larger and real-world instances}\label{EL}

\noindent {In this subsection, we present the results regarding the experiments with the larger subdivisions, IST500 and IST600, which have 500 and 600 curricula, and between 2500 to 2700 and 2800 to 3000 lectures, and with the real-world instance, with 1288 curricula, and 4102 lectures. }

\begin{figure}[h]
    \centering
    \includegraphics[scale=0.2]{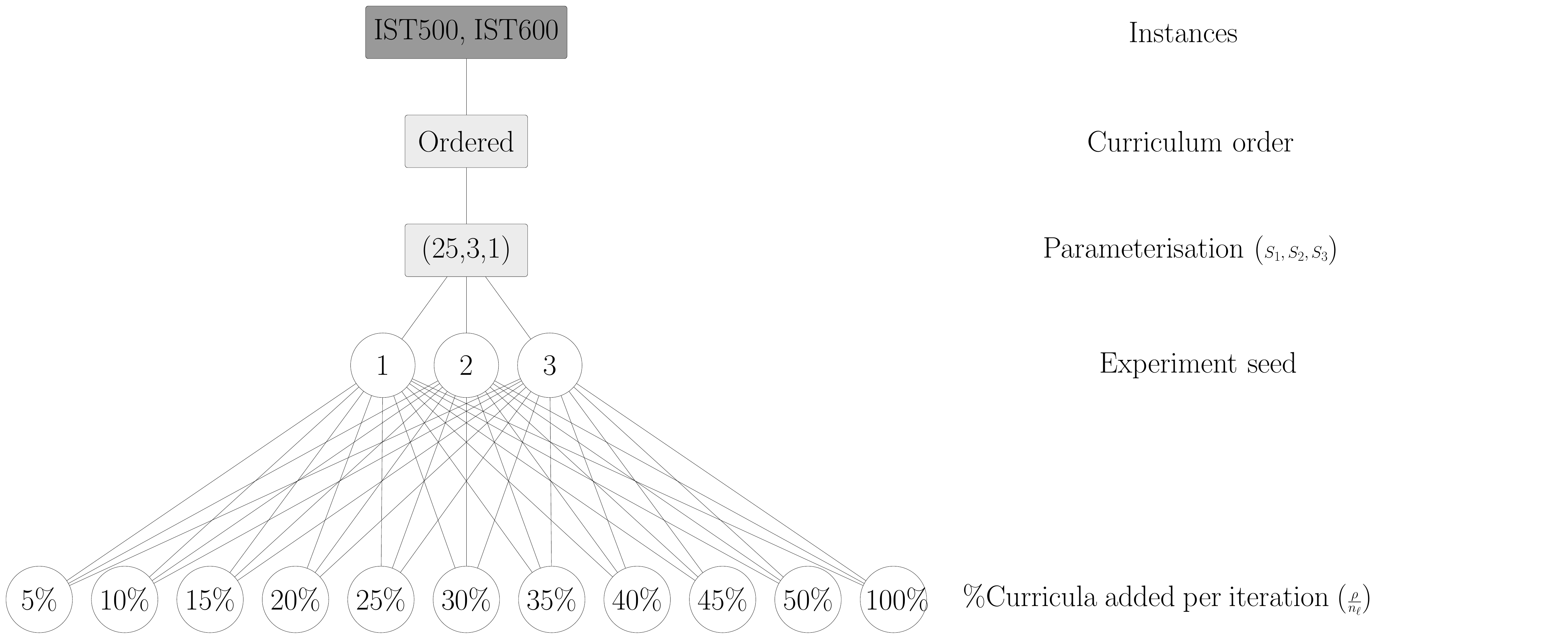}
    \caption{IST500 and IST600 experiments.}
    \label{exp_ls}
\end{figure}

\begin{figure}[h]
    \centering
    \includegraphics[scale=0.2]{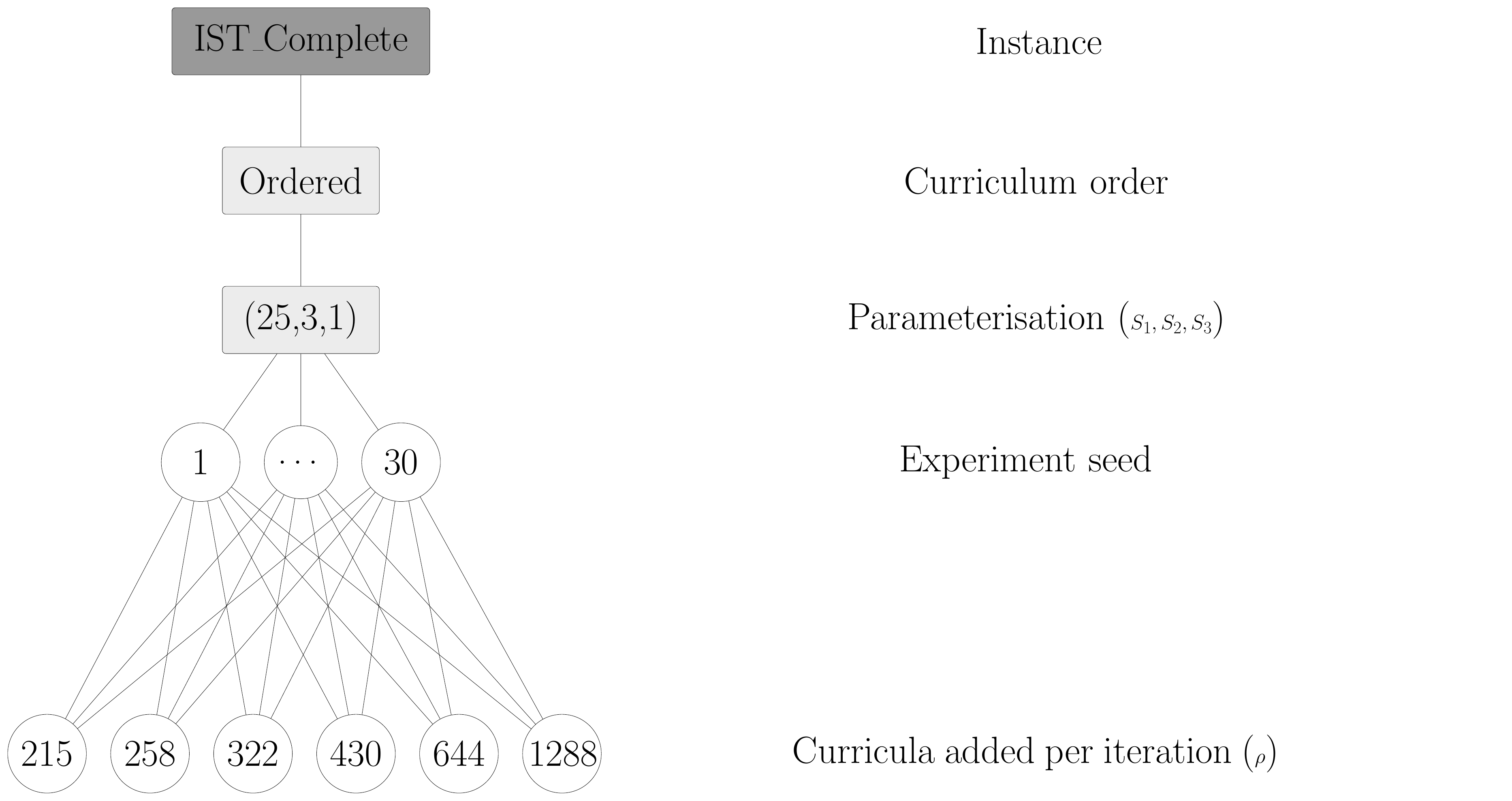}
    \caption{Real-world instance experiments.}
    \label{exp_rw}
\end{figure}

{There are some structural differences between the different groups of instances. For example, the IST400 subdivisions kept more than half of the lectures, while including less than a third of the original curricula. This indicates that in the original instance, and in larger subdivisions, the lectures are taught to more curricula on average. Since the subdivisions removed curricula and lectures, and not professors, the professors in the larger instances have considerably more lectures to teach. Finally, the larger subdivisions are increasingly more constrained regarding room availability, and the original instance also presents more room scarcity when compared to instances IST400. The experiments are summarised in Figure \ref{exp_ls} for the larger subdivisions and in Figure \ref{exp_rw} for the real-world instance. The experiments with 1288 curricula added per iteration are equivalent to solving the problem without decomposition.}

\begin{figure}[h]
    \centering
\begin{tikzpicture}[scale=0.7]
	\begin{axis}[%
  xmin=180,xmax=679,ymin=1200,ymax=1700,
  xtick distance = 100, 
  ytick distance = 100, 
 xlabel={Curricula added per iteration},
	ylabel={CPU(s)}]
\addplot coordinates {
	(215,	1631)
(258,	1412)
(322,	1465)
(430,	1501)
(644,	1276)

};

\addplot+[no marks] coordinates {
	(215,	1549)
(258,	1549)
(322,	1549)
(429,	1549)
(644,	1549)

};
\legend{{Decomposition}, Without decomposition}
	\end{axis}
\end{tikzpicture}
 \caption{Real-world experiment results.}
    \label{exp_results_rw}
\end{figure}
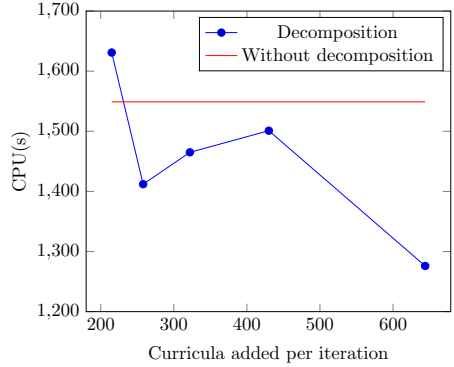

{Decomposition reduced the solution times by 12.2\% for instances IST500, by 8.9\% for instances IST600, and by 17.6\% for the real-world instance. For the experiments with instances IST400 and IST500, on average, the best curricula added per iteration values represented 20\% of the total curricula. However, for instances IST600 and the original instance, the best values are 25\% and 50\% of the total curricula being added per iteration. The results are presented in Figure \ref{exp_results_500} for instances IST500, in Figure \ref{exp_results_600} for instances IST600, and in Figure \ref{exp_results_rw} for the real-world instance.}

\begin{figure}[h]
\begin{minipage}{.5\textwidth}
    \centering
\begin{tikzpicture}[scale=0.7]
	\begin{axis}[%
  xmin=0,xmax=275,ymin=160,ymax=200,
  xtick distance = 50, 
  ytick distance = 10, 
 xlabel={Curricula added per iteration},
	ylabel={CPU(s)}]
\addplot coordinates {
	(25,	175)
(50,	169)
(75,	170)
(100,	165)
(125,	167)
(150,	166)
(175,	169)
(200,	166)
(225,	178)
(250,	177)

};

\addplot+[no marks] coordinates {
(25,	188)
(50,	188)
(75,	188)
(100,	188)
(125,	188)
(150,	188)
(175,	188)
(200,	188)
(225,	188)
(250,	188)

};

\legend{{Decomposition},{Without decomposition}}
	\end{axis}
\end{tikzpicture}
 \caption{IST500 experiment results.}
    \label{exp_results_500}
    \end{minipage}
\hspace{0.02cm}
\begin{minipage}{.5\textwidth}
    \centering
\begin{tikzpicture}[scale=0.7]
	\begin{axis}[%
  xmin=0,xmax=330,ymin=290,ymax=360,
  xtick distance = 60, 
  ytick distance = 10, 
 xlabel={Curricula added per iteration},
	ylabel={CPU(s)}]
\addplot coordinates {
	(30,	313)
(60,	321)
(90,	312)
(120,	311)
(150,	295)
(180,	329)
(210,	321)
(240,	324)
(270,	329)
(300,	305)

};

\addplot+[no marks] coordinates {
(30,	324)
(60,	324)
(90,	324)
(120,	324)
(150,	324)
(180,	324)
(210,	324)
(240,	324)
(270,	324)
(300,	324)

};

\legend{{Decomposition},{Without decomposition}}
	\end{axis}
\end{tikzpicture}
 \caption{IST600 experiment results.}
    \label{exp_results_600}
    \end{minipage}
\end{figure}

{The structural differences between the instances can explain this discrepancy in results. If a lecture is shared between multiple curricula, it is more likely that introducing curricula that include said lecture will create constraint violations in the decomposed problem. Likewise, if the professors have more lectures to teach, adding more lectures to the problem has a greater probability of making their timetables infeasible. In practice, this translates into larger increments being more efficient for the more complete instances, as they have more interdependencies between curricula.
Finally, in instances where there are fewer available rooms per lecture, the introduction of new lectures may violate the room availability constraints. All of these concerns translate into more changes being required to solve constraint violations when incrementing the problem, as each set of curricula is less contained in itself, \textit{i.e.}, there is more interaction between different curricula, reducing the relative effectiveness of the decomposition strategy.}

When solving a decomposed instance there is a trade-off between the amount of time required to adjust timetables which become infeasible after an increment, which is reduced by doing larger increments, and the simplicity of each iteration, which is increased by doing smaller increments. This simplicity promotes the possibility of using faster diversification rates and facilitates the identification of the time slots with the most violations and the most challenging constraint violation groups. In Table \ref{out_runs}, both the number of runs per $\rho$ value that took less than 700 seconds and more than 2300 seconds for the real-world instance is presented. In general, smaller increments have an increased likelihood of producing both very fast and very slow runs. We speculate that even for these parameterisations, sometimes the generated assignments do not require a significant amount of changes. In these cases, the added simplicity can promote very efficient runs. Nonetheless, the opposite can happen, and the generated assignments require several corrections, resulting in slower solution times. These results suggest that further improvements can be made to the methodology if tighter curricula clusters can be created, enabling smaller increments without requiring as much ineffective work as they would have fewer interactions amongst each other. 

\begin{table}[h]
    \centering
    \caption{Outlier runs per curricula added per iteration.}
   \begin{tabular}{|c|c|c|}
    \hline
        $\rho$ &  $<700$s&$>2300$s\\
        \hline
         215&3&4\\
         258&0&3\\
322&0&4\\
430&2&5\\
\textbf{644}&\textbf{0}&\textbf{2}\\
1288&0&2\\
\hline
    \end{tabular}

    \label{out_runs}
\end{table}

\subsection{Competition instances}\label{Benchmark}

\begin{figure}[h]
    \centering
    \includegraphics[scale=0.25]{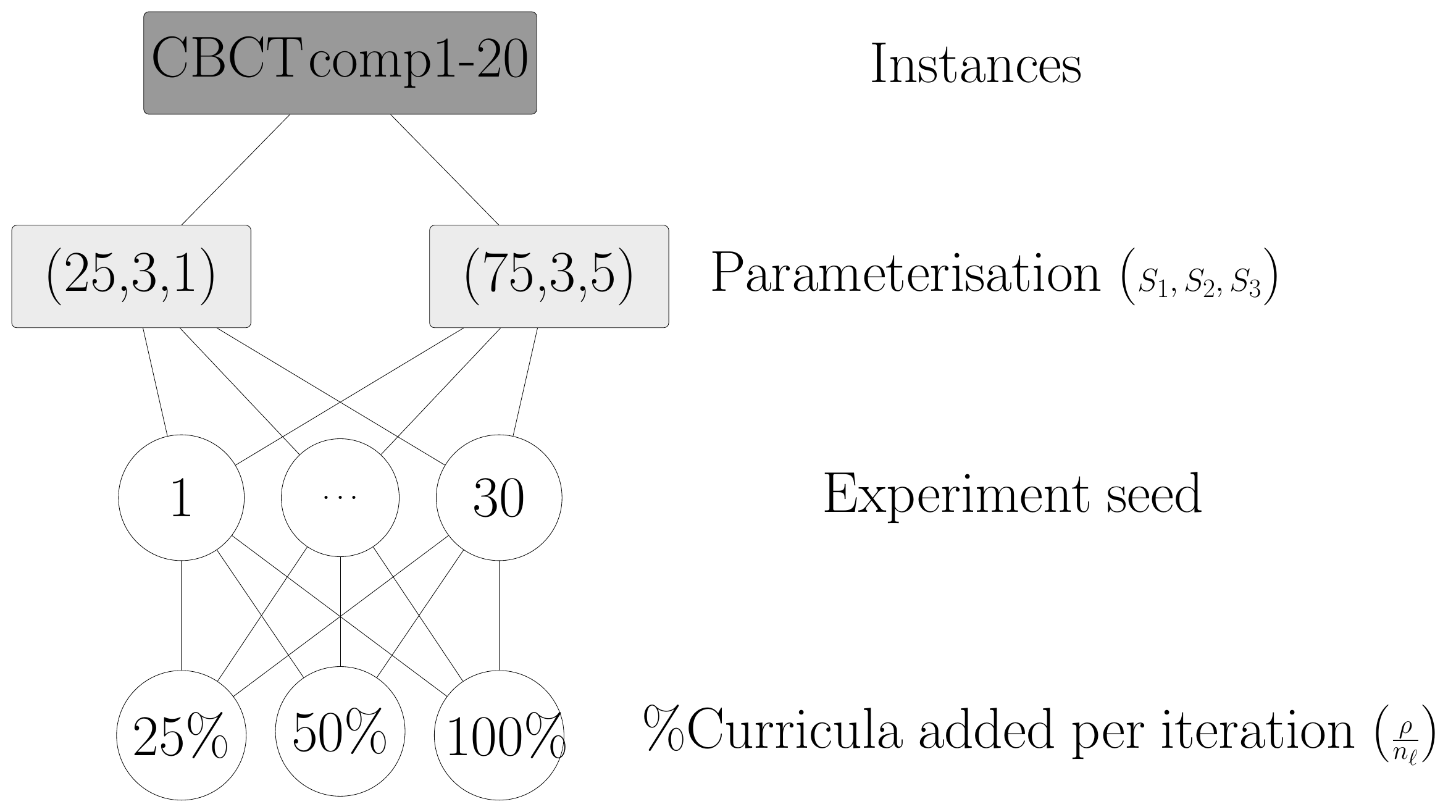}
    \caption{Adapted ITC 2007 CBCT instances experiments}
    \label{exp-itc}
\end{figure}

\noindent {In this subsection, we present the results regarding the experiments with the adapted benchmark instances from the CBCT competition track of ITC 2007. Maximum consecutive workload, precedence, and different day lecture constraints were added to the original instances, and unavailable slots were not regarded. Each of the twenty instances was run with 30 different seeds for parameterisations (25,3,1) and (75,3,5), with 25\%, 50\%, and 100\% of their curricula added per iteration. The experiments are summarised in Figure \ref{exp-itc}. }

The parameterisation (25,3,1) with 50\% of curricula being added per iteration was the fastest configuration for 15 out of the 20 tested instances. Since these instances are not very complex, this parameterisation, which diversifies the search faster, performed better. Moreover, making divisions smaller than 50\% of the curricula was also ineffective, as the number of curricula does not justify such small increments. 
The results are presented in Table \ref{avg bench}.

\begin{table}[h]
    \centering
     \caption{Average solution time for the adapted benchmark instances(s).}
    \label{avg bench}
   
    \begin{tabular}{|c|c|c|c|c|c|c|}
\hline
&\multicolumn{3}{c|}{(25,3,1)}& \multicolumn{3}{c|}{(75,3,5)}\\
\hline
Instance&100\%&50\%&25\%&100\%&50\%&25\%\\
\hline
comp1&0.044&\textbf{0.042}&0.053&0.048&0.048&0.051\\
comp2&1.283&\textbf{0.741}&0.926&1.429&2.197&0.761\\
comp3&0.442&0.419&\textbf{0.409}&0.463&0.509&0.443\\
comp4&\textbf{0.121}&0.137&0.137&0.130&0.132&0.134\\
comp5&3.089&\textbf{2.622}&2.841&3.866&3.693&4.099\\
comp6&0.416&0.384&\textbf{0.320}&0.548&0.349&0.390\\
comp7&0.674&\textbf{0.570}&0.601&0.780&0.734&0.598\\
comp8&0.169&\textbf{0.158}&0.189&0.168&0.182&0.202\\
comp9&0.426&\textbf{0.387}&0.421&0.527&0.494&0.403\\
comp10&0.461&\textbf{0.414}&0.473&0.594&0.848&0.445\\
comp11&0.103&\textbf{0.077}&0.082&0.090&0.101&0.090\\
comp12&4.067&\textbf{2.009}&2.444&2.784&2.508&2.688\\
comp13&0.470&\textbf{0.181}&0.211&0.204&0.190&0.224\\
comp14&0.491&\textbf{0.232}&0.234&0.246&0.235&0.243\\
comp15&0.538&0.441&\textbf{0.433}&0.467&0.510&0.481\\
comp16&0.411&\textbf{0.364}&0.405&0.615&0.475&0.433\\
comp17&0.435&0.381&\textbf{0.371}&0.579&0.457&0.389\\
comp18&0.090&\textbf{0.077}&0.090&0.092&0.082&0.096\\
comp19&0.447&\textbf{0.402}&0.433&0.420&0.417&0.459\\
comp20&0.570&\textbf{0.515}&\textbf{0.515}&0.672&0.630&0.593\\
\hline
    \end{tabular}
\end{table}

\section{Conclusion and future research}\label{c}

\noindent We propose a hybrid meta-heuristic for generating feasible course timetables. It uses elements from ALNS, GLS, and VNS. Moreover, a new instance decomposition and increment procedure is also introduced. The methodology is tested in real-world instances and significantly outperforms the commercial software currently used by the university.

Solving the decomposed problem promoted significant time reductions in both the original instance, 18\%, and its randomly generated subdivisions, 27\%. This methodology eases the identification of local optima by dealing with more optimised versions of the problem in each iteration. This allows faster diversification rates to be more efficient. Moreover, since our methodology focuses the search on the time slots with more constraint violations and attempts to determine which constraint groups are the hardest to reduce, having simpler problems facilitates these processes. 

However, decomposing the instance poses the drawback that some of the work done to find feasible schedules is fruitless, as later increments might reveal that the previous set of assignments is infeasible. Nonetheless, this effect can be reduced by clustering curricula that share some traits between themselves, such as professors or lectures, making the increments more contained. In the computational experiments, using ordered increments promoted solution time reductions 18\% greater than those obtained when incrementing the instances randomly. 

Regarding future research, several aspects of this work can be further explored. The methodology was applied to find feasible timetables. Nonetheless, it would be interesting to extend it to problems with objectives to be optimised. Moreover, applying the decomposition and increment procedure with different meta-heuristics and other problems should also be an interesting pursuit. Finally, the ordering algorithm implemented was simple, and it still significantly improved solution times. However, machine-learning techniques \citep{Zeitrag2022, Nagar2023} could be applied to find the best clustering for the decomposition and increment procedure.

\bmhead{Acknowledgments}

\noindent All the authors acknowledge the Portuguese national funds through the FCT - Foundation for Science and Technology, I.P., under the project UI/BD/154023/2022. Declarations of interest: none.


\bibliography{sn-bibliography}


\begin{thebibliography}{40}
\ifx \bisbn   \undefined \def \bisbn  #1{ISBN #1}\fi
\ifx \binits  \undefined \def \binits#1{#1}\fi
\ifx \bauthor  \undefined \def \bauthor#1{#1}\fi
\ifx \batitle  \undefined \def \batitle#1{#1}\fi
\ifx \bjtitle  \undefined \def \bjtitle#1{#1}\fi
\ifx \bvolume  \undefined \def \bvolume#1{\textbf{#1}}\fi
\ifx \byear  \undefined \def \byear#1{#1}\fi
\ifx \bissue  \undefined \def \bissue#1{#1}\fi
\ifx \bfpage  \undefined \def \bfpage#1{#1}\fi
\ifx \blpage  \undefined \def \blpage #1{#1}\fi
\ifx \burl  \undefined \def \burl#1{\textsf{#1}}\fi
\ifx \doiurl  \undefined \def \doiurl#1{\url{https://doi.org/#1}}\fi
\ifx \betal  \undefined \def \betal{\textit{et al.}}\fi
\ifx \binstitute  \undefined \def \binstitute#1{#1}\fi
\ifx \binstitutionaled  \undefined \def \binstitutionaled#1{#1}\fi
\ifx \bctitle  \undefined \def \bctitle#1{#1}\fi
\ifx \beditor  \undefined \def \beditor#1{#1}\fi
\ifx \bpublisher  \undefined \def \bpublisher#1{#1}\fi
\ifx \bbtitle  \undefined \def \bbtitle#1{#1}\fi
\ifx \bedition  \undefined \def \bedition#1{#1}\fi
\ifx \bseriesno  \undefined \def \bseriesno#1{#1}\fi
\ifx \blocation  \undefined \def \blocation#1{#1}\fi
\ifx \bsertitle  \undefined \def \bsertitle#1{#1}\fi
\ifx \bsnm \undefined \def \bsnm#1{#1}\fi
\ifx \bsuffix \undefined \def \bsuffix#1{#1}\fi
\ifx \bparticle \undefined \def \bparticle#1{#1}\fi
\ifx \barticle \undefined \def \barticle#1{#1}\fi
\bibcommenthead
\ifx \bconfdate \undefined \def \bconfdate #1{#1}\fi
\ifx \botherref \undefined \def \botherref #1{#1}\fi
\ifx \url \undefined \def \url#1{\textsf{#1}}\fi
\ifx \bchapter \undefined \def \bchapter#1{#1}\fi
\ifx \bbook \undefined \def \bbook#1{#1}\fi
\ifx \bcomment \undefined \def \bcomment#1{#1}\fi
\ifx \oauthor \undefined \def \oauthor#1{#1}\fi
\ifx \citeauthoryear \undefined \def \citeauthoryear#1{#1}\fi
\ifx \endbibitem  \undefined \def \endbibitem {}\fi
\ifx \bconflocation  \undefined \def \bconflocation#1{#1}\fi
\ifx \arxivurl  \undefined \def \arxivurl#1{\textsf{#1}}\fi
\csname PreBibitemsHook\endcsname

\bibitem[\protect\citeauthoryear{Tan et~al.}{2021}]{TanEtAl2021}
\begin{botherref}
\oauthor{\bsnm{Tan}, \binits{J.S.}},
\oauthor{\bsnm{Goh}, \binits{S.L.}},
\oauthor{\bsnm{Kendall}, \binits{G.}},
\oauthor{\bsnm{Sabar}, \binits{N.R.}}:
A survey of the state-of-the-art of optimisation methodologies in school timetabling problems.
Expert Systems with Applications
\textbf{165}
(2021)
\doiurl{10.1016/j.eswa.2020.113943}
\end{botherref}
\endbibitem

\bibitem[\protect\citeauthoryear{Leite et~al.}{2019}]{LEITE2019}
\begin{barticle}
\bauthor{\bsnm{Leite}, \binits{N.}},
\bauthor{\bsnm{Melício}, \binits{F.}},
\bauthor{\bsnm{C.~Rosa}, \binits{A.}}:
\batitle{A fast simulated annealing algorithm for the examination timetabling problem}.
\bjtitle{Expert Systems with Applications}
\bvolume{122},
\bfpage{137}--\blpage{151}
(\byear{2019})
\doiurl{10.1016/j.eswa.2018.12.048}
\end{barticle}
\endbibitem

\bibitem[\protect\citeauthoryear{Almeida et~al.}{2023}]{ALMEIDA2023}
\begin{barticle}
\bauthor{\bsnm{Almeida}, \binits{J.}},
\bauthor{\bsnm{Santos}, \binits{D.}},
\bauthor{\bsnm{Figueira}, \binits{J.R.}},
\bauthor{\bsnm{Francisco}, \binits{A.P.}}:
\batitle{A multi-objective mixed integer linear programming model for thesis defence scheduling}.
\bjtitle{European Journal of Operational Research}
(\byear{2023})
\doiurl{10.1016/j.ejor.2023.06.031}
\end{barticle}
\endbibitem

\bibitem[\protect\citeauthoryear{Babaei et~al.}{2015}]{BabaeiEtAl2015}
\begin{barticle}
\bauthor{\bsnm{Babaei}, \binits{H.}},
\bauthor{\bsnm{Karimpour}, \binits{J.}},
\bauthor{\bsnm{Hadidi}, \binits{A.}}:
\batitle{A survey of approaches for university course timetabling problem}.
\bjtitle{Computers and Industrial Engineering}
\bvolume{86},
\bfpage{43}--\blpage{59}
(\byear{2015})
\doiurl{10.1016/j.cie.2014.11.010}
\end{barticle}
\endbibitem

\bibitem[\protect\citeauthoryear{Thepphakorn and Pongcharoen}{2020}]{ThepphakornPongchaoren2020}
\begin{botherref}
\oauthor{\bsnm{Thepphakorn}, \binits{T.}},
\oauthor{\bsnm{Pongcharoen}, \binits{P.}}:
Performance improvement strategies on cuckoo search algorithms for solving the university course timetabling problem.
Expert Systems with Applications
\textbf{161}
(2020)
\doiurl{10.1016/j.eswa.2020.113732}
\end{botherref}
\endbibitem

\bibitem[\protect\citeauthoryear{Chen et~al.}{2021}]{ChenEtAl2021}
\begin{barticle}
\bauthor{\bsnm{Chen}, \binits{M.C.}},
\bauthor{\bsnm{Sze}, \binits{S.N.}},
\bauthor{\bsnm{Goh}, \binits{S.L.}},
\bauthor{\bsnm{Sabar}, \binits{N.R.}},
\bauthor{\bsnm{Kendall}, \binits{G.}}:
\batitle{A survey of university course timetabling problem: Perspectives, trends and opportunities}.
\bjtitle{IEEE Access}
\bvolume{9},
\bfpage{106515}--\blpage{106529}
(\byear{2021})
\doiurl{10.1109/ACCESS.2021.3100613}
\end{barticle}
\endbibitem

\bibitem[\protect\citeauthoryear{Vermuyten et~al.}{2016}]{VermuytenEtAl2016}
\begin{barticle}
\bauthor{\bsnm{Vermuyten}, \binits{H.}},
\bauthor{\bsnm{Lemmens}, \binits{S.}},
\bauthor{\bsnm{Marques}, \binits{I.}},
\bauthor{\bsnm{Beliën}, \binits{J.}}:
\batitle{Developing compact course timetables with optimized student flows}.
\bjtitle{European Journal of Operational Research}
\bvolume{251}(\bissue{2}),
\bfpage{651}--\blpage{661}
(\byear{2016})
\doiurl{10.1016/j.ejor.2015.11.028}
\end{barticle}
\endbibitem

\bibitem[\protect\citeauthoryear{Boufflet et~al.}{2021}]{BouffletEtAl2021}
\begin{botherref}
\oauthor{\bsnm{Boufflet}, \binits{J.-P.}},
\oauthor{\bsnm{Arbaoui}, \binits{T.}},
\oauthor{\bsnm{Moukrim}, \binits{A.}}:
The student scheduling problem at université de technologie de compiègne.
Expert Systems with Applications
\textbf{175}
(2021)
\doiurl{10.1016/j.eswa.2021.114735}
\end{botherref}
\endbibitem

\bibitem[\protect\citeauthoryear{Herres and Schmitz}{2021}]{Herres2021}
\begin{barticle}
\bauthor{\bsnm{Herres}, \binits{B.}},
\bauthor{\bsnm{Schmitz}, \binits{H.}}:
\batitle{Decomposition of university course timetabling: A systematic study of subproblems and their complexities}.
\bjtitle{Annals of Operations Research}
\bvolume{302}(\bissue{2}),
\bfpage{405}--\blpage{423}
(\byear{2021})
\doiurl{10.1007/s10479-019-03382-0}
\end{barticle}
\endbibitem

\bibitem[\protect\citeauthoryear{Soria-Alcaraz et~al.}{2016}]{Soria-Alcaraz2016}
\begin{barticle}
\bauthor{\bsnm{Soria-Alcaraz}, \binits{J.A.}},
\bauthor{\bsnm{Özcan}, \binits{E.}},
\bauthor{\bsnm{Swan}, \binits{J.}},
\bauthor{\bsnm{Kendall}, \binits{G.}},
\bauthor{\bsnm{Carpio}, \binits{M.}}:
\batitle{Iterated local search using an add and delete hyper-heuristic for university course timetabling}.
\bjtitle{Applied Soft Computing Journal}
\bvolume{40},
\bfpage{581}--\blpage{593}
(\byear{2016})
\doiurl{10.1016/j.asoc.2015.11.043}
\end{barticle}
\endbibitem

\bibitem[\protect\citeauthoryear{Kiefer et~al.}{2017}]{Kiefer2017}
\begin{barticle}
\bauthor{\bsnm{Kiefer}, \binits{A.}},
\bauthor{\bsnm{Hartl}, \binits{R.F.}},
\bauthor{\bsnm{Schnell}, \binits{A.}}:
\batitle{Adaptive large neighborhood search for the curriculum-based course timetabling problem}.
\bjtitle{Annals of Operations Research}
\bvolume{252}(\bissue{2}),
\bfpage{255}--\blpage{282}
(\byear{2017})
\doiurl{10.1007/s10479-016-2151-2}
\end{barticle}
\endbibitem

\bibitem[\protect\citeauthoryear{Bagger et~al.}{2019a}]{BaggerEtAl2019B}
\begin{barticle}
\bauthor{\bsnm{Bagger}, \binits{N.-C.F.}},
\bauthor{\bsnm{Sørensen}, \binits{M.}},
\bauthor{\bsnm{Stidsen}, \binits{T.R.}}:
\batitle{Dantzig–wolfe decomposition of the daily course pattern formulation for curriculum-based course timetabling}.
\bjtitle{European Journal of Operational Research}
\bvolume{272}(\bissue{2}),
\bfpage{430}--\blpage{446}
(\byear{2019})
\doiurl{10.1016/j.ejor.2018.06.042}
\end{barticle}
\endbibitem

\bibitem[\protect\citeauthoryear{Bagger et~al.}{2019b}]{BaggerEtAl2019}
\begin{barticle}
\bauthor{\bsnm{Bagger}, \binits{N.-C.F.}},
\bauthor{\bsnm{Kristiansen}, \binits{S.}},
\bauthor{\bsnm{Sørensen}, \binits{M.}},
\bauthor{\bsnm{Stidsen}, \binits{T.R.}}:
\batitle{Flow formulations for curriculum-based course timetabling}.
\bjtitle{Annals of Operations Research}
\bvolume{280}(\bissue{1-2}),
\bfpage{121}--\blpage{150}
(\byear{2019})
\doiurl{10.1007/s10479-018-3096-4}
\end{barticle}
\endbibitem

\bibitem[\protect\citeauthoryear{Pillay and Özcan}{2019}]{PillayOzcan2019}
\begin{barticle}
\bauthor{\bsnm{Pillay}, \binits{N.}},
\bauthor{\bsnm{Özcan}, \binits{E.}}:
\batitle{Automated generation of constructive ordering heuristics for educational timetabling}.
\bjtitle{Annals of Operations Research}
\bvolume{275}(\bissue{1}),
\bfpage{181}--\blpage{208}
(\byear{2019})
\doiurl{10.1007/s10479-017-2625-x}
\end{barticle}
\endbibitem

\bibitem[\protect\citeauthoryear{Lindahl et~al.}{2019}]{LindahlEtAl2019}
\begin{barticle}
\bauthor{\bsnm{Lindahl}, \binits{M.}},
\bauthor{\bsnm{Stidsen}, \binits{T.}},
\bauthor{\bsnm{Sørensen}, \binits{M.}}:
\batitle{Quality recovering of university timetables}.
\bjtitle{European Journal of Operational Research}
\bvolume{276}(\bissue{2}),
\bfpage{422}--\blpage{435}
(\byear{2019})
\doiurl{10.1016/j.ejor.2019.01.026}
\end{barticle}
\endbibitem

\bibitem[\protect\citeauthoryear{Gülcü and Akkan}{2020}]{GülcüAkkan2020}
\begin{barticle}
\bauthor{\bsnm{Gülcü}, \binits{A.}},
\bauthor{\bsnm{Akkan}, \binits{C.}}:
\batitle{Robust university course timetabling problem subject to single and multiple disruptions}.
\bjtitle{European Journal of Operational Research}
\bvolume{283}(\bissue{2}),
\bfpage{630}--\blpage{646}
(\byear{2020})
\doiurl{10.1016/j.ejor.2019.11.024}
\end{barticle}
\endbibitem

\bibitem[\protect\citeauthoryear{Akkan et~al.}{2021}]{AkkanEtAl2021}
\begin{botherref}
\oauthor{\bsnm{Akkan}, \binits{C.}},
\oauthor{\bsnm{Gülcü}, \binits{A.}},
\oauthor{\bsnm{Kuş}, \binits{Z.}}:
Minimum penalty perturbation heuristics for curriculum-based timetables subject to multiple disruptions.
Computers and Operations Research
\textbf{132}
(2021)
\doiurl{10.1016/j.cor.2021.105306}
\end{botherref}
\endbibitem

\bibitem[\protect\citeauthoryear{de~la Rosa-Rivera et~al.}{2021}]{Rosa-Rivera2021}
\begin{botherref}
\oauthor{\bsnm{Rosa-Rivera}, \binits{F.}},
\oauthor{\bsnm{Nunez-Varela}, \binits{J.I.}},
\oauthor{\bsnm{Ortiz-Bayliss}, \binits{J.C.}},
\oauthor{\bsnm{Terashima-Marín}, \binits{H.}}:
Algorithm selection for solving educational timetabling problems.
Expert Systems with Applications
\textbf{174}
(2021)
\doiurl{10.1016/j.eswa.2021.114694}
\end{botherref}
\endbibitem

\bibitem[\protect\citeauthoryear{Song et~al.}{2021}]{Song2021}
\begin{botherref}
\oauthor{\bsnm{Song}, \binits{T.}},
\oauthor{\bsnm{Chen}, \binits{M.}},
\oauthor{\bsnm{Xu}, \binits{Y.}},
\oauthor{\bsnm{Wang}, \binits{D.}},
\oauthor{\bsnm{Song}, \binits{X.}},
\oauthor{\bsnm{Tang}, \binits{X.}}:
Competition-guided multi-neighborhood local search algorithm for the university course timetabling problem.
Applied Soft Computing
\textbf{110}
(2021)
\doiurl{10.1016/j.asoc.2021.107624}
\end{botherref}
\endbibitem

\bibitem[\protect\citeauthoryear{Goh et~al.}{2017}]{GohEtAl2017}
\begin{barticle}
\bauthor{\bsnm{Goh}, \binits{S.L.}},
\bauthor{\bsnm{Kendall}, \binits{G.}},
\bauthor{\bsnm{Sabar}, \binits{N.R.}}:
\batitle{Improved local search approaches to solve the post enrolment course timetabling problem}.
\bjtitle{European Journal of Operational Research}
\bvolume{261}(\bissue{1}),
\bfpage{17}--\blpage{29}
(\byear{2017})
\doiurl{10.1016/j.ejor.2017.01.040}
\end{barticle}
\endbibitem

\bibitem[\protect\citeauthoryear{Nagata}{2018}]{Nagata2018}
\begin{barticle}
\bauthor{\bsnm{Nagata}, \binits{Y.}}:
\batitle{Random partial neighborhood search for the post-enrollment course timetabling problem}.
\bjtitle{Computers and Operations Research}
\bvolume{90},
\bfpage{84}--\blpage{96}
(\byear{2018})
\doiurl{10.1016/j.cor.2017.09.014}
\end{barticle}
\endbibitem

\bibitem[\protect\citeauthoryear{Goh et~al.}{2019}]{GohEtAl2019}
\begin{barticle}
\bauthor{\bsnm{Goh}, \binits{S.L.}},
\bauthor{\bsnm{Kendall}, \binits{G.}},
\bauthor{\bsnm{Sabar}, \binits{N.R.}}:
\batitle{Simulated annealing with improved reheating and learning for the post enrolment course timetabling problem}.
\bjtitle{Journal of the Operational Research Society}
\bvolume{70}(\bissue{6}),
\bfpage{873}--\blpage{888}
(\byear{2019})
\doiurl{10.1080/01605682.2018.1468862}
\end{barticle}
\endbibitem

\bibitem[\protect\citeauthoryear{Sylejmani et~al.}{2022}]{Sylejmani2022}
\begin{barticle}
\bauthor{\bsnm{Sylejmani}, \binits{K.}},
\bauthor{\bsnm{Gashi}, \binits{E.}},
\bauthor{\bsnm{Ymeri}, \binits{A.}}:
\batitle{Simulated annealing with penalization for university course timetabling}.
\bjtitle{Journal of Scheduling}
(\byear{2022})
\doiurl{10.1007/s10951-022-00747-5}
\end{barticle}
\endbibitem

\bibitem[\protect\citeauthoryear{Bettinelli et~al.}{2015}]{BettinelliEtAl2015}
\begin{barticle}
\bauthor{\bsnm{Bettinelli}, \binits{A.}},
\bauthor{\bsnm{Cacchiani}, \binits{V.}},
\bauthor{\bsnm{Roberti}, \binits{R.}},
\bauthor{\bsnm{Toth}, \binits{P.}}:
\batitle{An overview of curriculum-based course timetabling}.
\bjtitle{TOP}
\bvolume{23}(\bissue{2}),
\bfpage{313}--\blpage{349}
(\byear{2015})
\doiurl{10.1007/s11750-015-0366-z}
\end{barticle}
\endbibitem

\bibitem[\protect\citeauthoryear{Akkan and Gülcü}{2018}]{AkkanGulcu2018}
\begin{barticle}
\bauthor{\bsnm{Akkan}, \binits{C.}},
\bauthor{\bsnm{Gülcü}, \binits{A.}}:
\batitle{A bi-criteria hybrid genetic algorithm with robustness objective for the course timetabling problem}.
\bjtitle{Computers and Operations Research}
\bvolume{90},
\bfpage{22}--\blpage{32}
(\byear{2018})
\doiurl{10.1016/j.cor.2017.09.007}
\end{barticle}
\endbibitem

\bibitem[\protect\citeauthoryear{Wu}{2011}]{Wu2011}
\begin{barticle}
\bauthor{\bsnm{Wu}, \binits{C.-C.}}:
\batitle{Parallelizing a clips-based course timetabling expert system}.
\bjtitle{Expert Systems with Applications}
\bvolume{38}(\bissue{6}),
\bfpage{7517}--\blpage{7525}
(\byear{2011})
\doiurl{10.1016/j.eswa.2010.12.116}
\end{barticle}
\endbibitem

\bibitem[\protect\citeauthoryear{Shiau}{2011}]{Shiau2011}
\begin{barticle}
\bauthor{\bsnm{Shiau}, \binits{D.-F.}}:
\batitle{A hybrid particle swarm optimization for a university course scheduling problem with flexible preferences}.
\bjtitle{Expert Systems with Applications}
\bvolume{38}(\bissue{1}),
\bfpage{235}--\blpage{248}
(\byear{2011})
\doiurl{10.1016/j.eswa.2010.06.051}
\end{barticle}
\endbibitem

\bibitem[\protect\citeauthoryear{Aladag et~al.}{2009}]{Aladag2009}
\begin{barticle}
\bauthor{\bsnm{Aladag}, \binits{C.H.}},
\bauthor{\bsnm{Hocaoglu}, \binits{G.}},
\bauthor{\bsnm{Basaran}, \binits{M.A.}}:
\batitle{The effect of neighborhood structures on tabu search algorithm in solving course timetabling problem}.
\bjtitle{Expert Systems with Applications}
\bvolume{36}(\bissue{10}),
\bfpage{12349}--\blpage{12356}
(\byear{2009})
\doiurl{10.1016/j.eswa.2009.04.051}
\end{barticle}
\endbibitem

\bibitem[\protect\citeauthoryear{Bagger et~al.}{2018}]{Bagger2018EtAl}
\begin{barticle}
\bauthor{\bsnm{Bagger}, \binits{N.-C.F.}},
\bauthor{\bsnm{Sørensen}, \binits{M.}},
\bauthor{\bsnm{Stidsen}, \binits{T.R.}}:
\batitle{Benders’ decomposition for curriculum-based course timetabling}.
\bjtitle{Computers and Operations Research}
\bvolume{91},
\bfpage{178}--\blpage{189}
(\byear{2018})
\doiurl{10.1016/j.cor.2017.10.009}
\end{barticle}
\endbibitem

\bibitem[\protect\citeauthoryear{Lindahl et~al.}{2018}]{LindahlEtAl2018}
\begin{barticle}
\bauthor{\bsnm{Lindahl}, \binits{M.}},
\bauthor{\bsnm{Mason}, \binits{A.J.}},
\bauthor{\bsnm{Stidsen}, \binits{T.}},
\bauthor{\bsnm{Sørensen}, \binits{M.}}:
\batitle{A strategic view of university timetabling}.
\bjtitle{European Journal of Operational Research}
\bvolume{266}(\bissue{1}),
\bfpage{35}--\blpage{45}
(\byear{2018})
\doiurl{10.1016/j.ejor.2017.09.022}
\end{barticle}
\endbibitem

\bibitem[\protect\citeauthoryear{Dunke and Nickel}{2023}]{Dunke2023}
\begin{barticle}
\bauthor{\bsnm{Dunke}, \binits{F.}},
\bauthor{\bsnm{Nickel}, \binits{S.}}:
\batitle{A matheuristic for customized multi-level multi-criteria university timetabling}.
\bjtitle{Annals of Operations Research}
\bvolume{328}(\bissue{2}),
\bfpage{1313}--\blpage{1348}
(\byear{2023})
\doiurl{10.1007/s10479-023-05325-2}
\end{barticle}
\endbibitem

\bibitem[\protect\citeauthoryear{Sousa et~al.}{2016}]{Sousa2016}
\begin{barticle}
\bauthor{\bsnm{Sousa}, \binits{T.}},
\bauthor{\bsnm{Morais}, \binits{H.}},
\bauthor{\bsnm{Castro}, \binits{R.}},
\bauthor{\bsnm{Vale}, \binits{Z.}}:
\batitle{Evaluation of different initial solution algorithms to be used in the heuristics optimization to solve the energy resource scheduling in smart grids}.
\bjtitle{Applied Soft Computing Journal}
\bvolume{48},
\bfpage{491}--\blpage{506}
(\byear{2016})
\doiurl{10.1016/j.asoc.2016.07.028}
\end{barticle}
\endbibitem

\bibitem[\protect\citeauthoryear{Juman and Hoque}{2015}]{Juman2015}
\begin{barticle}
\bauthor{\bsnm{Juman}, \binits{Z.A.M.S.}},
\bauthor{\bsnm{Hoque}, \binits{M.A.}}:
\batitle{An efficient heuristic to obtain a better initial feasible solution to the transportation problem}.
\bjtitle{Applied Soft Computing Journal}
\bvolume{34},
\bfpage{813}--\blpage{826}
(\byear{2015})
\doiurl{10.1016/j.asoc.2015.05.009}
\end{barticle}
\endbibitem

\bibitem[\protect\citeauthoryear{Amaliah et~al.}{2022}]{Amaliah2022}
\begin{botherref}
\oauthor{\bsnm{Amaliah}, \binits{B.}},
\oauthor{\bsnm{Fatichah}, \binits{C.}},
\oauthor{\bsnm{Suryani}, \binits{E.}}:
A supply selection method for better feasible solution of balanced transportation problem.
Expert Systems with Applications
\textbf{203}
(2022)
\doiurl{10.1016/j.eswa.2022.117399}
\end{botherref}
\endbibitem

\bibitem[\protect\citeauthoryear{Viana et~al.}{2022}]{Viana2022}
\begin{botherref}
\oauthor{\bsnm{Viana}, \binits{B.M.F.}},
\oauthor{\bsnm{Pereira}, \binits{L.T.}},
\oauthor{\bsnm{Toledo}, \binits{C.F.M.}},
\oauthor{\bsnm{Santos}, \binits{S.R.}},
\oauthor{\bsnm{Maia}, \binits{S.M.D.M.}}:
Feasible–infeasible two-population genetic algorithm to evolve dungeon levels with dependencies in barrier mechanics.
Applied Soft Computing
\textbf{119}
(2022)
\doiurl{10.1016/j.asoc.2022.108586}
\end{botherref}
\endbibitem

\bibitem[\protect\citeauthoryear{Song et~al.}{2018}]{Song2018}
\begin{barticle}
\bauthor{\bsnm{Song}, \binits{T.}},
\bauthor{\bsnm{Liu}, \binits{S.}},
\bauthor{\bsnm{Tang}, \binits{X.}},
\bauthor{\bsnm{Peng}, \binits{X.}},
\bauthor{\bsnm{Chen}, \binits{M.}}:
\batitle{An iterated local search algorithm for the university course timetabling problem}.
\bjtitle{Applied Soft Computing Journal}
\bvolume{68},
\bfpage{597}--\blpage{608}
(\byear{2018})
\doiurl{10.1016/j.asoc.2018.04.034}
\end{barticle}
\endbibitem

\bibitem[\protect\citeauthoryear{Pillay and Özcan}{2019}]{Pillay2019}
\begin{barticle}
\bauthor{\bsnm{Pillay}, \binits{N.}},
\bauthor{\bsnm{Özcan}, \binits{E.}}:
\batitle{Automated generation of constructive ordering heuristics for educational timetabling}.
\bjtitle{Annals of Operations Research}
\bvolume{275}(\bissue{1}),
\bfpage{181}--\blpage{208}
(\byear{2019})
\doiurl{10.1007/s10479-017-2625-x}
\end{barticle}
\endbibitem

\bibitem[\protect\citeauthoryear{Goh et~al.}{2019}]{Goh2019}
\begin{barticle}
\bauthor{\bsnm{Goh}, \binits{S.L.}},
\bauthor{\bsnm{Kendall}, \binits{G.}},
\bauthor{\bsnm{Sabar}, \binits{N.R.}}:
\batitle{Monte carlo tree search in finding feasible solutions for course timetabling problem}.
\bjtitle{International Journal on Advanced Science, Engineering and Information Technology}
\bvolume{9}(\bissue{6}),
\bfpage{1936}--\blpage{1943}
(\byear{2019})
\doiurl{10.18517/ijaseit.9.6.10224}
\end{barticle}
\endbibitem

\bibitem[\protect\citeauthoryear{Zeitrag et~al.}{2022}]{Zeitrag2022}
\begin{botherref}
\oauthor{\bsnm{Zeitrag}, \binits{Y.}},
\oauthor{\bsnm{Figueira}, \binits{J.R.}},
\oauthor{\bsnm{Horta}, \binits{N.}},
\oauthor{\bsnm{Neves}, \binits{R.}}:
Surrogate-assisted automatic evolving of dispatching rules for multi-objective dynamic job shop scheduling using genetic programming.
Expert Systems with Applications
\textbf{209}
(2022)
\doiurl{10.1016/j.eswa.2022.118194}
\end{botherref}
\endbibitem

\bibitem[\protect\citeauthoryear{Nagar et~al.}{2023}]{Nagar2023}
\begin{botherref}
\oauthor{\bsnm{Nagar}, \binits{J.}},
\oauthor{\bsnm{Chaturvedi}, \binits{S.K.}},
\oauthor{\bsnm{Soh}, \binits{S.}},
\oauthor{\bsnm{Singh}, \binits{A.}}:
A machine learning approach to predict the k-coverage probability of wireless multihop networks considering boundary and shadowing effects.
Expert Systems with Applications
\textbf{226}
(2023)
\doiurl{10.1016/j.eswa.2023.120160}
\end{botherref}
\endbibitem

\end{thebibliography}

\begin{appendices}
\section{Increment algorithms}

\begin{algorithm}[h]
\caption{Fixed curriculum increment: \textit{CurriculumIncrement()}}\label{alg_fixed_inc}

	\begin{algorithmic}[1]\normalsize
	\State{\textbf{input:} $s^*$, $\rho$, $L$, $L'$,$I'$};

        \If{($ length(L') + \rho \leqslant length(L)$)}{\;$\ell'\leftarrow \rho$;}
        \Else{\; $\ell'\leftarrow length(L)-length(L')$};
        \EndIf
        \For{($\ell=1,\ldots,\ell'$)}\State{$\ell''\leftarrow L[length(L')]$;}
        \State{$L'\leftarrow L'\cup\ell''$;}
        \State{$s^*\leftarrow s^*\cup ConstructInitialTimetable(\ell'', s^*,I')$;}
        \EndFor
        
        \end{algorithmic}
\end{algorithm}

\begin{algorithm}[h]
\caption{Violations based curriculum increment: \textit{CurriculumIncrement()}}\label{alg_vio_inc}
	\begin{algorithmic}[1]\normalsize
	\State{\textbf{input:} $s^*$,  $\overline{f}$, $L$, $L'$, $I'$};
	
        \While{($length(L') < length(L) \land f(s^*)<\overline{f}$)}
        \State{$\ell''\leftarrow L[length(L')]$;}
        \State{$L'\leftarrow L'\cup \ell''$;}
        \State{$s^*\leftarrow s^*\cup ConstructInitialTimetable(\ell'', s^*,I')$;}
        \EndWhile
        \end{algorithmic} 
\end{algorithm}

\newpage

\section{Initial solutions}

\begin{algorithm}[h]
\caption{\textit{ConstructInitialTimetable()}}\label{alg_initial}
\begin{algorithmic}[1]\normalsize
	\State{\textbf{input:} $ct$, $\ell$, $\overline{j}$, $\overline{k}$, $w'_j$, $r'_{ijk}$, $I'$};
	\State{\textbf{output:} $ct$};
        \State{$j\leftarrow rand(1,\ldots,\overline{j})$};
        \State{$k\leftarrow rand(1,\ldots,\overline{k})$};
   
        \For{($i \in a_{\ell} \land i \notin I' $)}
        \State{$assigning \leftarrow true$;}
       
        \While{($assigning$)}
        
        \While{($w'_j+c_i>m_\ell \lor r'_{ijk}>r_\ell$)}{\;$k\leftarrow k+1$};
        \If{($j>n_j \lor k+c_i>n_k $)}
        \State{$j\leftarrow j + 1$};
        \If{($j>n_j$)}{\;$j\leftarrow 1$};
        \EndIf
        \State{$k\leftarrow rand(1,\ldots,\overline{k})$};
        \EndIf
        \EndWhile

        \State{$assigning \leftarrow false$;}

        \For{($k' \in \{k,\ldots,k + c_i-1\}$)}

        \If{($y_{jk'\ell}=1$)}{\;$assigning\leftarrow true$};
        \EndIf
        \EndFor
        
        \EndWhile
       
        \State{$x_{ijk}\leftarrow 1$};
        \State{$ct \leftarrow ct \cup x_{ijk}$};
        \State{$I'\leftarrow I' \cup i$};
        \State{$w'_j\leftarrow w'_j + c_i$};
        \State{$r'_{i',j,k+c_i}\leftarrow r'_{i',j,k+c_i} + r'_{ijk} + c_i$};
        \State{$k\leftarrow k+c_i$};
        \If{($k + c_{i+1}>n_k$)}
        \State{$j\leftarrow j + 1$};
        \If{($j>n_j$)}{\;$j\leftarrow 1$};
        \EndIf
        \State{$k\leftarrow rand(1,\ldots,\overline{k})$};
        \EndIf

        \EndFor

        \end{algorithmic} 
	   
\end{algorithm}

\newpage

\section{Neighbourhood structures}

\begin{algorithm}[h]
\caption{Worst slot neighbourhood search}\label{alg_worst_slot}
\begin{algorithmic}[1]\normalsize
	\State{\textbf{input:} $s^*$, $swap$, $I''$};
	\State{\textbf{output:} $s'$};

        \State{$i'\leftarrow rand(I'')$};
        \State{$s'\leftarrow s^*$};
        
        \If{($swap$)}
        \State{$\ell'\leftarrow rand(L'')$};
        \State{$i''\leftarrow rand(a_{\ell'})$};
        \State{$x'_{i'j''k''}\leftarrow 1$};
        \State{$x'_{i''j'k'}\leftarrow 1$};
        \Else

        \State{$j \leftarrow rand(1,\ldots,n_j)$};
        \State{$k \leftarrow rand(1,\ldots,n_k-c_{i'})$};
        \State{$x'_{i'jk}\leftarrow 1$};
        \EndIf
        \State{{\textbf{return}($\mbox{$s'$}$);}}
        \end{algorithmic} 
	   
\end{algorithm}

\begin{algorithm}[h]
\caption{Worst curriculum neighbourhood search}\label{alg_worst_curr}
	\begin{algorithmic}[1]\normalsize
	\State{\textbf{input:} $s^*$, $swap$, $L''$};
	\State{\textbf{output:} $s'$};

        \State{$\ell'\leftarrow rand(L'')$};

        \State{$\hat{j}\leftarrow MostViolationsDay(\ell')$};
        \State{$\hat{k}\leftarrow MostViolationsHour(\ell')$};
        \State{$I''\leftarrow \{\}$};
        
        \For{($i''\in a_{\ell'}$)}
        \For{($k\in  \{\hat{k}-c_{i''}+1,\ldots,\hat{k}\}$)}
        \If{($x_{i''\hat{j}k}=1$)}{\;$I''\leftarrow I''\cup i''$};
        \EndIf
        \EndFor
        \EndFor
        
        \State{$i'\leftarrow rand(I'')$};
        \State{$s'\leftarrow s^*$};
        
        \If{($swap$)}
        \State{$i''\leftarrow rand(a_{\ell'})$};
        \State{$x'_{i'j''k''}\leftarrow 1$};
        \State{$x'_{i''j'k'}\leftarrow 1$};
        \Else

        \State{$j \leftarrow rand(1,\ldots,n_j)$};
        \State{$k \leftarrow rand(1,\ldots,n_k-c_i)$};
        \State{$x'_{i'jk}\leftarrow 1$};
        \EndIf
        \State{{\textbf{return}($\mbox{$s'$}$);}}
        \end{algorithmic} 
	   
\end{algorithm}

\begin{algorithm}[h]
\caption{Worst professor neighbourhood search}\label{alg_worst_prof}
\begin{algorithmic}[1]\normalsize
	\State{\textbf{input:} $s^*$, $swap$, $P''$, $L$};
	\State{\textbf{output:} $s'$};

        \State{$p'\leftarrow rand(P'')$};

        \State{$\hat{j}\leftarrow MostViolationsDay(p')$};
        \State{$\hat{k}\leftarrow MostViolationsHour(p')$};
        \State{$I''\leftarrow \{\}$};
        
        \For{($i''\in b_{p'}$)}
        \For{($k\in  \{\hat{k}-c_{i''}+1,\ldots,\hat{k}\}$)}
        \If{($x_{i''\hat{j}k}=1$)}{\;$I''\leftarrow I''\cup i''$};
        \EndIf
        \EndFor
        \EndFor
        
        \State{$i'\leftarrow rand(I'')$};
        \State{$s'\leftarrow s^*$};
        
        \If{($swap$)}
        \State{$\ell'\leftarrow rand(L, i'\in a_{\ell'})$};
        \State{$i''\leftarrow rand(a_{\ell'})$};
        \State{$x'_{i'j''k''}\leftarrow 1$};
        \State{$x'_{i''j'k'}\leftarrow 1$};
        \Else

        \State{$j \leftarrow rand(1,\ldots,n_j)$};
        \State{$k \leftarrow rand(1,\ldots,n_k-c_i)$};
        \State{$x'_{i'jk}\leftarrow 1$};
        \EndIf
        \State{{\textbf{return}($\mbox{$s'$}$);}}
        \end{algorithmic} 
\end{algorithm}

\begin{algorithm}[h]
\caption{Worst room type neighbourhood search}\label{alg_worst_room}
\begin{algorithmic}[1]\normalsize
	\State{\textbf{input:} $s^*$, $swap$, $\mu$, $t$, $\hat{j}$, $\hat{k}$, $I'$};
	\State{\textbf{output:} $s'$};

        \State{$I''\leftarrow \{\}$};
        
        \For{($i''\in g_\mu \land v_{i''t}=1$)}
        \For{($k\in  \{\hat{k}-c_{i''}+1,\ldots,\hat{k}\}$)}
        \If{($x_{i''\hat{j}k}=1$)}{\;$I''\leftarrow I''\cup i''$};
        \EndIf
        \EndFor
        \EndFor
        
        \State{$i'\leftarrow rand(I'')$};
        \State{$s'\leftarrow s^*$};
        
        \If{($swap$)}
        \State{$\ell'\leftarrow rand(L, i'\in a_{\ell'})$};
        \State{$i''\leftarrow rand(a_{\ell'}, v_{i''t}=0)$};
        \State{$x'_{i'j''k''}\leftarrow 1$};
        \State{$x'_{i''j'k'}\leftarrow 1$};
        \Else

        \State{$j \leftarrow rand(1,\ldots,n_j)$};
        \State{$k \leftarrow rand(1,\ldots,n_k-c_{i'})$};
        \State{$x'_{i'jk}\leftarrow 1$};
        \EndIf
        \State{{\textbf{return}($\mbox{$s'$}$);}}
        \end{algorithmic} 
	   
\end{algorithm}

\begin{algorithm}[h]
\caption{Different day neighbourhood search}\label{alg_diff_da}
\begin{algorithmic}[1]\normalsize
	\State{\textbf{input:} $s^*$, $swap$, $i'$};
	\State{\textbf{output:} $s'$};

        \State{$s'\leftarrow s^*$};
        
        \If{($swap$)}
        \State{$\ell'\leftarrow rand(L'' \land i'\in a_{\ell'})$};
        \State{$i''\leftarrow rand(a_{\ell'} \land \sum_{k=1}^{n_k} x_{i''j'k})=0$};
        \State{$x'_{i'j''k''}\leftarrow 1$};
        \State{$x'_{i''j'k'}\leftarrow 1$};
        \Else

        \State{$j \leftarrow rand(1,\ldots,n_j \land j\neq j')$};
        \State{$k \leftarrow rand(1,\ldots,n_k-c_{i'})$};
        \State{$x'_{i'jk}\leftarrow 1$};
        \EndIf
        \State{{\textbf{return}($\mbox{$s'$}$);}}
        \end{algorithmic} 
\end{algorithm}

\begin{algorithm}[h]
\caption{Precedence neighbourhood search}\label{alg_prec}
	\begin{algorithmic}[1]\normalsize
	\State{\textbf{input:} $s^*$, $swap$, $i'$, $i''$};
	\State{\textbf{output:} $s'$};

        \State{$s'\leftarrow s^*$};
        
        \If{($swap$)}
       
        \State{$x'_{i'j''k''}\leftarrow 1$};
        \State{$x'_{i''j'k'}\leftarrow 1$};
        \Else

        \State{$j \leftarrow rand(1,\ldots,n_j \land j\leqslant j')$};
        \If{($j<j'$)}
            \State{$k \leftarrow rand(1,\ldots,n_k-c_{i'})$};
        \Else
             \State{$k \leftarrow rand(1,\ldots,k'')$};
        \EndIf
        \State{$x'_{i'jk}\leftarrow 1$};
        \EndIf
        \State{{\textbf{return}($\mbox{$s'$}$);}}
        \end{algorithmic} 
	   
\end{algorithm}




\end{appendices}


\end{document}